\documentclass{gOMS2e}

\usepackage{cmap}					
\usepackage{mathtext}

\usepackage{epstopdf}
\usepackage{subfigure}

\theoremstyle{plain}
\newtheorem{theorem}{Theorem}[section]

\newtheorem{lemma}[theorem]{Lemma}

\theoremstyle{definition}

\theoremstyle{remark}

\usepackage{amsmath,amsfonts,amssymb,amsthm} 
\usepackage{icomma}

\usepackage{graphicx}  
\graphicspath{{images/}{images2/}}  
\usepackage{wrapfig} 

\usepackage{array,tabularx,tabulary,booktabs} 
\usepackage{longtable}  
\usepackage{multirow} 

\usepackage{etoolbox}

\usepackage{lastpage} 

\usepackage{soul}

\usepackage{multicol}

\usepackage[linesnumbered,ruled,vlined]{algorithm2e}
\usepackage[usenames]{color}

\newcommand{\bR}{\mathbb{R}}

\DeclareMathOperator*{\argmin}{argmin}
\DeclareMathOperator*{\argmax}{argmax}

\newtheorem{prop}{Proposition}

{
\theoremstyle{remark}

}
{
\theoremstyle{definition}

}

\usepackage{hyperref}
\usepackage{verbatim}
\usepackage{ragged2e}
\allowdisplaybreaks

\begin{document}


\articletype{Optimization}

\title{Efficient Numerical Methods to Solve Sparse Linear Equations with Application to PageRank}
 
\author{
\name{Anton Anikin\textsuperscript{a}\thanks{ Anton Anikin Email: anikin@icc.ru }, 
Alexander Gasnikov\textsuperscript{b,c,d}\thanks{ Alexander Gasnikov. Email: gasnikov@yandex.ru},
Alexander Gornov\textsuperscript{a}\thanks{ Alexander Gornov Email: gornov@icc.ru}, 
Dmitry Kamzolov\textsuperscript{b}\thanks{ Dmitry Kamzolov Email: kamzolov.dmitry@phystech.edu}, \\
Yury Maximov \textsuperscript{e,f}\thanks{ Yury Maximov Email: yury@lanl.gov} and
Yurii Nesterov \textsuperscript{g}\thanks{ Yurii Nesterov Email: nesterov@core.ucl.ac.be}}
\affil{\textsuperscript{a} Institute of System Dynamics and Control Theory (ISDCT SB RAS), Russia;\\
\textsuperscript{b}Moscow Institute of Physics and Technology, Moscow, Russia;\\
\textsuperscript{c}Higher School of Economics, Moscow, Russia;\\
\textsuperscript{d}Weierstrass Institute for Applied Analysis and Stochastics, Berlin, Germany;\\
\textsuperscript{e} Skolkovo Institute of Science and Technology, Russia;\\
\textsuperscript{f} Theoretical Division T-5, Los Alamos National Laboratory, USA;\\
\textsuperscript{g}Center
 for Operations Research and Econometrics (CORE),\\ Catholic University of Louvain (UCL), Belgium.
}
\received{v1.0 released May 2019}
}

\maketitle
\begin{abstract}
Over the last two decades, the PageRank problem has received increased interest from the academic community as an efficient tool to estimate web-page importance in information retrieval. 
Despite numerous developments, the design of efficient optimization algorithms for the PageRank problem is still a challenge. This paper proposes three new algorithms with a linear-time complexity for solving the problem over a bounded-degree graph. The idea behind them is to set up the PageRank as a convex minimization problem over a unit simplex, and then solve it using iterative methods with small iteration complexity. Our theoretical results are supported by an extensive empirical justification using real-world and simulated data.  
\end{abstract}


\begin{classcode}
90C25, 90C47, 90C60
\end{classcode}

\section{Introduction}
\label{intro}
In this paper, we aim at solving a system of linear equations $Px = b$ for a sufficiently sparse matrix $A$ and a vector $b$, $P \in \mathbb{R}^{m\times n}$ and $b\in \mathbb{R}^m$,  which is primarily motivated by finding a stationary distribution of a Markov Chain over a sparse large-scale communication graph. This problem, known as the PageRank, was pioneered by Brin and Page in \cite{brin2012reprint,page1999pagerank} in the early nineties; however, it still attracts significant interest from both the academic community and the industry.

In the PageRank we assume, that the (asymmetric) transition probability matrix $P$, associated with the web-graph, is known; so, the problem is to find a stationary distribution or a vector $x$, $x\in \mathbb{R}^n$, with the coordinates corresponding to an expected portion of time which a random walk spends at a particular node, $x = P^\top x$. 

The problem becomes especially challenging in high dimensions, since direct computations of an inverse matrix become inefficient due to non-linear time and memory efforts. Many studies have been devoted to the approximation of the PageRank vector based on a random walk analysis and Markov Chain Monte Carlo methods~\cite{avrachenkov2007monte,gasnikov2015efficient,kamvar2004adaptive,kamvar2003extrapolation,sarma2015fast,tong2006fast}. Those methods are very attractive, both theoretically and practically, while the spectral gap, i.e., the difference between the two largest eigenvalues of the transition matrix, is sufficiently large. The latter is often not the case for sparse graphs with complex topologies, and, furthermore, estimating the gap requires significant time and effort, as well~\cite{levin2017markov}. 

Another line of research aims at finding a  stationary distribution of a Markov chain using convex optimization~\cite{gasnikov2015efficiency,nazin2011randomized,nesterov2014subgradient,nesterov2015finding}. According to these results, the PageRank can be equally stated as an ~$\ell_p$-norm minimization, $p \ge 1$, constrained on a unit simplex: 
\begin{gather}\label{eq:main}
\|P^\top x - x\|_p \to \min_{x\in \Delta_1^n}, 
\end{gather}
where $\Delta_1^n = \{x: \sum_{i=1}^n x^i = 1, x^i \ge 0\}$. 

From the convex optimization perspective, similar problems appear in applied mathematics, statistics, and machine learning. Among these are LASSO \cite{friedman2001elements}, a traffic matrix estimation in large-scale communication networks~\cite{zhang2003fast}, phase recovery in a linearized model of electric current~\cite{stott2009dc}, and a finite element method~\cite{logan2011first}.
The high-dimensional nature of the problems above call into a question of the utility of traditional approaches, which do not devote sufficient attention to the problem structure and require  non-linear time and memory efforts. 

In this paper, we focus on the influence of the transition probability matrix sparsity on the computational complexity of the PageRank problem. We advocate particular efficiency of convex optimization methods for the simplex constrained $\|Ax\|_2^2$ and $\|Ax\|_\infty$  minimization problems~\cite{bubeck2015convex,candes2008enhancing}. In particular, we prove that the time complexity of these problems is linear in the problem dimension if the number of non-zeros in each row/column is bounded above by some constant $d$. 
The key contribution is a set of efficient algorithms to update a function value, a gradient, and an argument in an (almost) dimension-independent manner. 
Later in the paper, we extend our results to a more general setup, where a limited number of dense rows or columns in the transition matrix is allowed. 

\subsection{Contribution}
Our contribution is as follows. We propose: 
\begin{enumerate}
\item The \texttt{NL1} algorithm, a $\ell_1$-proximal gradient descent method that supports sparse updates of the gradient and the function value. It allows us to solve the problem with overall time 
complexity $O(nd^2 \log^2 n + d^2\log^2 (n/\varepsilon)/\varepsilon^2)$, where $d$ is the maximal in- and out- vertexes' degree of the graph, 
$\varepsilon$ is required accuracy, and $\|Px - x\|_2 \le \varepsilon$;
\item The \texttt{S-FW} method, an extension of the Frank--Wolfe algorithm that allows efficient updates to the gradient and the objective value. The resulting algorithm \texttt{S-FW} has better time-complexity estimates compared to the \texttt{NL1} method, $O(n + d^2 \log (2+n/d^2)/\varepsilon^2)$. Also, the algorithm often has remarkably better performance in practice;
\item The \texttt{GK} algorithm, which aims to minimize $\|P^\top x - x\|_{\infty}$ over the unit simplex. We provide an equivalent saddle-point setup for this problem and, subsequently, solve it subsequently  by a version of 
the mirror descent with a randomized projecting. We prove that a randomized projection with KL-divergence guarantees the total running time of the algorithm to be bounded from above by $O(n + d \log n \log(n /\delta)/\varepsilon^2)$ with probability at least $1 - \delta$, for any $0 < \delta < 1$;
\item Finally, we extend the linear-time complexity estimates to the sparse graphs with a small number of dense  rows/columns containing more than $d$ non-zeros. 
\end{enumerate}

Let us emphasize that a $d$-sparse matrix could have as many non-zero elements as $n\cdot d$. Time-complexity estimates of the proposed algorithms are yet sub-linear in the number of non-zero elements of the transition matrix if its sparsity pattern is known beforehand. The last statement implies that there is no need to read all non-zero elements of a stochastic matrix to arrive at approximate solution of the PageRank problem. In Table~\ref{table:page-rank-time} we summarize the best known results for solving the PageRank problem using optimization techniques.

\begin{table}[th]\footnotesize
  \centering
  \begin{tabular}{| l | c | p{3.8cm} | p{2.2cm} |}
      \hline
      {\bf Algorithm} & {\bf Constraint} & {\bf Time Complexity} & {\bf Objective} \\
      \hline
      Nazin - Polyak \cite{nazin2011randomized} & no & $O\left(\frac{n \log (n/\delta)}{\varepsilon^2}\right)$ & $\|P^\top x - x\|_2 \le \varepsilon$ \hspace{-4mm}\\
      \hline
      Nesterov \cite{nesterov2012coordinate} & $d$-sparse & $O\left(n + \frac{d^2 \log n}{\varepsilon^2}{}\right)$ & $\mathbb{E} \|P^\top x - x\|_2 \le \varepsilon$ \hspace{-4mm} \\
      \hline
      Nesterov \cite{nesterov2014subgradient} & $d$-sparse & $O\left(dn \frac{\log n}{\varepsilon^2}\right)$& $\|P^\top x - x\|_\infty \le \varepsilon$\hspace{-4mm}\\ 
      && $O\left(\frac{d^{1/2}n^{3/2} \log n}{\varepsilon}\right)$ & \hspace{-4mm}\\
      \hline
      Juditsky et al. \cite{juditsky2011first,nemirovski2009robust}& no & $O\left(\frac{n \log (n/\delta)}{\varepsilon^2}\right)$ & $\|P^\top x - x\|_\infty \le \varepsilon$ \hspace{-4mm}\\
      \hline
      Nesterov - Nemirovski \cite{nesterov2015finding} \hspace{-2mm} & $d$-sparse & $O\left(\frac{nd}{\alpha} \log 1/\varepsilon \right)$ & $\|x - x^*\|_1 \le \varepsilon$ \hspace{-4mm}\\
      & Google page & & \\
      \hline
      Polyak - Tremba \cite{polyak2012regularization} & $d$-sparse & $O\left(\frac{dn}{\varepsilon}\right)$ & $\|P^\top x - x\|_1 \le \varepsilon$\hspace{-4mm}\\
      & on average && \\
      \hline 
      Gasnikov - Dmitriev, \cite{gasnikov2015efficient} & $d$-sparse & $O\left(n + \frac{d\log n \log(n/\delta)}{\beta \varepsilon^2}\right)$\hspace{-1mm} & $\mathbb{E}\|x - x^*\|_2 \le \varepsilon$ \hspace{-4mm}\\
      & \hspace{-3.5mm}spectral gap $\beta$\hspace{-4mm} &&\\
      \hline
      Gasnikov - Dmitriev, \cite{gasnikov2015efficient} & $m$ non-zeros & $O\left(m + \left(n + \frac{m^2}{n^2}\!\right)\!\frac{\log n}{\varepsilon^2}\!\!\right)$\hspace{-1mm} & $\|P^T x - x\|_\infty \le \varepsilon$ \hspace{-4mm}\\
      & &&\\
      \hline
Langville - Meyer, \cite{langville2011google} & $d$-sparse & $O\left(\frac{dn}{\beta} \log \frac{n}{\varepsilon}\right)$ & $\|x - x^*\|_1 \le \varepsilon$ \hspace{-4mm}\\
      & \hspace{-3.5mm}spectral gap $\beta$\hspace{-4mm} &&
      \\
      \hline
Cohen et al. \cite{cohen2017almost} & $m$ non-zeros & $O(m + n^{1+o(1)}) \log^{O(1)}  \left(\frac{n}{\beta\varepsilon}\right)$ & $\mathbb{E} \|P^\top x - x\|_2 \le \varepsilon$ \hspace{-4mm}\\
      & &&
      \\ 
      \hline
      This paper, Section~\ref{sec:NL1}  & $d$-sparse & $O\left(n d^2 \log^{2} n + \frac{d^2 \log(n/d^2 + 1)}{\varepsilon^2}\right)$& $\|P^\top x - x\|_2 \le \varepsilon$ \hspace{-4mm}\\ 
      &  &&\\
\hline
      This paper, Section~\ref{sec:s-fw} & $d$-sparse & $O\left(n  +  \frac{d^2 \log (n/d^2 + 2)}{\varepsilon^2}\right)$\hspace{-2.5mm} & $\|P^\top x - x\|_2 \le \varepsilon$ \hspace{-4mm}\\ \hline
This paper, Section~\ref{sec:GK}  & $d$-sparse & $O\left(n  +  \frac{d \log n \log(n/\delta)}{\varepsilon^2}\right)$\hspace{-2.5mm} & $\|P^\top x - x\|_\infty \le \varepsilon$ \hspace{-4mm}\\
      \hline
  \end{tabular}
  \caption{Time complexity of the PageRank problem. Time complexity of the algorithms proposed in this paper along with results of Nesterov \cite{nesterov2012coordinate} are the only sub-linear algorithms known to the authors. Algorithm \cite{cohen2017almost} can be favourable in theory for high dimensions, but useless in practice due to a high degree of the logarithm.}
  \label{table:page-rank-time}
\end{table}

The Google page condition, required in \cite{nesterov2015finding}, implies the existence of a column~$j$ such that all $P_{ij} \ge \alpha$, $\alpha > 0$ for any $i:\, 1\le i \le n$; and $\beta$ corresponds to the spectral gap of the matrix $P$, e.g. $\beta = \lambda_1 (P) - \lambda_2(P) = 1 - \lambda_2(P)$, the difference between the largest and the second-largest eigenvalues of the transition matrix~$P$. Algorithm~\cite{nesterov2012coordinate} runs on average for the randomized coordinate descent algorithm while the estimates \cite{gasnikov2015efficient,nazin2011randomized,nemirovski2009robust} are correct with probability at least $1 - \delta$, for any~$\delta:\, 0 < \delta < 1$.

\subsection{Paper structure and notation} In Section~\ref{sec:NL1} we introduce a $\ell_1$ proximal gradient descent algorithm. The key idea behind this time-efficient algorithm, referenced below as \texttt{NL1}, is a sparse update to a gradient and a function value updates, $f(x) = \|Ax\|_2^2$, $A = P^\top - I$:
\[\nabla f(x+h) = \nabla f(x) + A^\top A h, \quad f(x+h) = f(x) + h^\top \nabla f(x) + \|Ah\|_2^2.\]
The \texttt{NL1} method does not require computation of a (full) gradient on each step if matrix $A$ is sufficiently sparse. We show that the update vector $h$ on each step has only two non-zero coordinates, which correspond to the minimal and maximal components of the gradient. To update gradient coordinates and extract the minimal and the maximal value, we use a list of binary heaps \cite{cormen2001introduction}, which admits logarithmic dependence of iteration complexity in dimension.

Algorithm \texttt{S-FW}, a revision of the Frank--Wolfe conditional gradient, is proposed in Section~\ref{sec:s-fw}. The Frank--Wolfe algorithm has recently stimulated much interest, mainly due to the numerous Big Data problems to which it has been applied \cite{harchaoui2015conditional,jaggi2013revisiting,nesterov2015complexity}. In this paper, we focus on an efficient gradient and a function value update for each iteration, which reduces the algorithm's overall time complexity. Accurate theoretical analysis results in a better time complexity estimate compared to~the~\texttt{NL1} algorithm. Also, the \texttt{S-FW} and \texttt{NL1} algorithms appear to be comparable from in practice. 

In Section~\ref{sec:GK} we deal with the problem of the approximation of the PageRank vector in the $\ell_\infty$ norm, i.e., the minimization of $\|P^\top x - x\|_\infty$. Using a technique of Juditsky et al. \cite{juditsky2011first}, the minimization can be equally stated as a saddle-point problem over the product of two unit simplexes. 
A randomized mirror descent seems to be one of the most popular tools for solving the problem efficiently. Its application 
to the PageRank improves time-complexity of an iteration to~$O(n \log (n/\delta) / \varepsilon^2)$, and does not depend on the problem sparsity~\cite{nemirovski2009robust}. 

Finally, in this paper, we investigate an approach pioneered by Grigoriadis and Khachiyan in \cite{grigoriadis1995sublinear}. The approach's idea is a randomized projection of the gradient on the simplex instead of a randomized approximation of the gradient itself. The gradient's projection on the unit simplex is carried out with the Kullback-Leibler divergence, corresponding to the exponential weighting of the gradient coordinates. The (sparse) randomized projection chooses one of the vertices of the unit simplex in such a way that the expected value gives an unbiased estimation of the projection itself~\cite{gasnikov2015efficient}. 

To this end, we design an algorithm to update the gradient vector on each iteration in an (almost) dimension-independent manner, e.g., with an (almost) linear time-complexity
\[O(n + d \log n \log (n/\sigma)/\varepsilon^2)\]
for $d$-sparse transition matrices, with $\|P^\top x- x\|_\infty \le \varepsilon$. Finally, we pay a special attention to its interpretation in terms of game theory~\cite{myerson2013game}. 

In Sections~\ref{sec:NL1} - \ref{sec:GK}, we summarize sparsification techniques allowing to reduce ranking problems over a dense graph to the problem on a sparse graph. The proposed algorithms deliver state-of-the-art time-complexity estimates for $d$-sparse optimization problems if $d$ is sufficiently small. Faster methods are possible via the Markov Chain Monte Carlo method for several situations. 
We refer to recent results~\cite{gasnikov2015efficient} for details. 

We conclude in Section~\ref{sec:case} with the implementation details and a case study. Finally, some technical proofs are given in the appendix.

In this paper, we use the following notation. By $\|x\|_p$, we denote the $\ell_p$
norm of vector $x \in \mathbb{R}^n$, in particular, $\|x\|_2 = x^\top x$  denotes the Euclidean norm, $\|x\|_1  = \sum_{i=1}^n |x_i|$, and $\|x\|_\infty= \max\{|x_1|, \dots , |x_n|\}$. 
By $O(1)$, we denote a positive constant. 

\section{\texorpdfstring{$\ell_1$}{L1}-Gradient descent for d-sparse problems}
\label{sec:NL1}
The PageRank problem can be equally stated as the following convex optimization problem for a sufficiently large positive constant $\gamma$
\begin{gather}\label{eq:pgm-setup}
f_\gamma(x) = \frac{1}{2} \|Ax\|_2^2 + \sum_{i = 1}^n\frac{\gamma}{2} (-x^i)_+^2 \to \min_{x^\top e = 1}, \qquad \gamma > 0,
\end{gather}
where $A = I - P^\top$, $I$ is an $n\times n$ identity matrix, $P \in \bR^{n\times n}_+$ is a column stochastic transition matrix, $e$ is an $n$-dimensional vector with each coordinate equals to one. Notice, that $f_\gamma(x)$ is a convex, non-negative, and  increasing function with $f_\gamma(x_*) = 0$, and $z_+ = \max\{z,0\}$. 

To solve the problem, we use a proximal gradient descent with an $\ell_1$ setup:
\begin{gather}\label{eq:pgm-step}
x_{k+1} = x_k + \argmin_{h: h^\top e = 0} \left\{f_\gamma(x_k) + h^\top \nabla f_\gamma(x_k) + \frac{L_1}{2} \|h\|_1^2\right\},
\end{gather}
where $f_\gamma(x)$ is $L_1$-smooth in $\ell_1$ norm
\[
	\|\nabla f_\gamma(x) - \nabla f_\gamma(y)\|_\infty \le L_1 \|x-y\|_1.
\]
In the case of the PageRank, $L_1 = 1 + \gamma$, as $P$ is a stochastic matrix. 
Lemma \ref{lem:nl1-min} implies that the vector $h_k$, which solves the problem \eqref{eq:pgm-step}, is always sparse whatever the function $f_\gamma(x)$ is. 

\begin{lemma}\label{lem:nl1-min}
A set of solutions of the following minimization problem
\begin{gather}\label{eq:l1-norm-min}
    \phi(h) = f_\gamma(x_k) + h^\top \nabla f_\gamma(x_k) + \frac{L_1}{2} \|h\|_1^2 \to \min_{h: h^\top e = 0}
\end{gather}
contains at least one vector $h_k$ with only two non-zero coordinates 
\[
h_k^{i_+} = - \frac{1}{4L_1}\left(\frac{\partial f_\gamma(x_k)}{\partial x^{i_+}} - \frac{\partial f_\gamma(x_k)}{\partial x^{i_-}}\right)  \text{ and } h_k^{i_-} =  \frac{1}{4L_1}\left(\frac{\partial f_\gamma(x_k)}{\partial x^{i_+}} - \frac{\partial f_\gamma(x_k)}{\partial x^{i_-}}\right),
\]
where 
$i_+ = \argmax_{1\le i\le n} \partial f_\gamma(x_k)/\partial x^i$ and $i_- = \argmin_{1\le i\le n} \partial f_\gamma(x_k)/\partial x^i$. 
\end{lemma}
\begin{proof}
Let $\bar h$ be a solution to the problem~\ref{eq:l1-norm-min}, so that $\|\bar h\|_1 = c$. The minimum to the linear function $f_\gamma(x_k) + h^\top \nabla f_\gamma(x_k) + \frac{L_1}{2} c^2 \to \min_{h: h^\top e = 0, \|h\|_1 = c}$ is attained at $h_k$ which has exactly two non-zero coordinates $h_k^{i^+} = - c/2$, $h_k^{i^-} = c/2$. The objective value is then \[\phi(h_k) = f_\gamma(x_k) - (\max_{1\le i\le n} \partial f_\gamma(x_k)/\partial x^i - \min_{1\le i\le n} \partial f_\gamma(x_k)/\partial x^i)c/2 + L_1 c^2/2.\]
Taking minimum in $c \ge 0$ we get the statement of the lemma. 
\end{proof}


By Lemma~\ref{lem:nl1-min}, the gradient $\nabla f_\gamma(x_{k+1})$ is 
\begin{align}\label{eq:pgm-grad-upd}
\nabla f_\gamma(x_{k+1}) = & \nabla f_\gamma(x_k + h_k)  \nonumber\\ 
= & A^\top A x_k + A^\top A h_k + \gamma \sum_{i = i_+, i_-}\left\{(-x_k^{i} - h_k^{i})_+^2 - (-x_k^i)_+^2\right\}.
\end{align}

The update vector, $h_k$, 
has at most $O(\min(n, d^2))$ non-zero coordinates. 
To efficiently update the gradient, we will use a doubly-linked list of binary heaps, so that the $j$-th heap is used to extract the minimal value and update coordinates from $(j-1)\lfloor n/d^2\rfloor + 1$ to $\min\{j\lfloor n/d^2\rfloor, n\}$ inclusively. We refer to a binary heap~\cite{cormen2001introduction} as a \textit{Max-Heap} (resp. {\it Min-Heap}), if the key stored in each node is greater (resp. less) than or equal to the keys in the node's children. We require both the minimal and the maximal coordinates of the gradient for an iteration of the algorithm. Using \textit{Min-Heap} (\textit{Max-Heap}) extracting the minimal (maximal) element from a heap of $m$ items requires $O(\log m)$ time. An update of a single element to preserve the keys' order requires $O(\log m)$ time, as well~\cite{cormen2001introduction}. That is why a gradient update~\eqref{eq:pgm-grad-upd} requires~$O(d^2 \log (2 + n/d^2))$ time in total.  

\begin{algorithm}[!th]
\DontPrintSemicolon 
\KwIn{$d$-sparse transition matrix $P$, starting point $x_0$,\; 
objective $f_\gamma(x) = \frac{1}{2} \|Ax\|_2^2 + \frac{\gamma}{2} \sum_{i=1}^n (-x_i)_+^2$ with $A = I - P^\top$ and $\gamma > 0$,\; 
number of iterations ${N}$ required by Theorem~\ref{thm:NL1} and accuracy $\varepsilon$. 
}
\KwOut{$x_{N}$, $\|Ax_N\|_2^2 \le \varepsilon^2$ where $N$ is given by Theorem \ref{thm:NL1}}
  \While{$k \le {N}$ and $f_\gamma(x) > \varepsilon^2$} {
     $i_+ = \argmax_{1\le i\le n} {\partial f_\gamma(x_k)}/{\partial x^i}$, \quad $i_- = \argmin_{1\le i\le n} {\partial f_\gamma(x_k)}/{\partial x^i}$\;
     $h_{k}^i = \begin{cases}
     			\frac{1}{4L_1}\left(\frac{\partial f_\gamma(x_k)}{\partial x^{i_+}} - \frac{\partial f_\gamma(x_k)}{\partial x^{i_-}}\right), & i = i_+\\
                - \frac{1}{4L_1}\left(\frac{\partial f_\gamma(x_k)}{\partial x^{i_+}} - \frac{\partial f_\gamma(x_k)}{\partial x^{i_-}}\right), & i = i_-\\
                0, & \text{otherwise}
     		   \end{cases}$\;
               Update argument: $x_{{k}+1} = x_{k} + {h_{k}}$\;
     Update gradient: $\nabla f_\gamma(x_{{k}+1}) = \nabla f_\gamma(x_{k}) + A^\top A h_{k} + \gamma (-x_k + h_k)_+ - \gamma(-x_k)_+$\;
     Update $f_\gamma(x_{{k}+1})$: \, $f_\gamma(x_{{k}+1}) = f_\gamma(x_{k}) + {h_{k}}^\top A^\top A {x_{k}} + \|A{h_{k}}\|_2^2/2 + \frac{\gamma}{2}\sum_{i: \, {h_{k}}^i \neq 0} (-{x_{k}}^i)_+^2 - (-x^i_{k} - {h_{k}}^i_{k})_+^2$\;
     ${k} = {k} + 1$\;
  }
\Return{$x_{k}$}\;
\caption{{\sc \texttt{NL1}: $\ell_1$ Gradient Descent for  PageRank}}
\label{algo:l1-pgm}
\end{algorithm}

Algorithm~\ref{algo:l1-pgm} presents the \texttt{NL1} algorithm with convergence rate established in Theorem~\ref{thm:NL1}.

\begin{theorem}\label{thm:NL1}
For a starting point $x_0$ in one of the vertices of the unit simplex, Algorithm~\ref{algo:l1-pgm} converges to $f_\gamma(x) \le \varepsilon^2$ for any constant $\gamma > 0$, $x \ge 0$, $e^\top x = 1$, with the overall time complexity
\[
T  = O\left(
    n d^2 \log (n/d^2 + 2){\log n} + \frac{d^2 \log (n/d^2 + 2)}{\varepsilon^2}
    \right).
\]
\end{theorem}

The complete proof of the theorem in given in the Appendix. 
The bound established in Theorem~\ref{thm:NL1} gives 
a sub-linear time complexity for the case $n \varepsilon^2 = o(1)$, 
while it is less practical for sufficiently large $\varepsilon$, $n \varepsilon^2 = \Omega(1)$. 

\textit{Matrix Sparsification.} In a number of practical problems, both column and row sparsity of the transition probability matrix seems to be very restrictive. Indeed, for the PageRank problem, search engines such as \texttt{Google} or \texttt{Yahoo!} may refer to a large number of sites simultaneously. In particular, that means that matrix $P^\top$ may have a few dense columns. Below we propose a method to sparsify the transition probability matrix $P$ in order to improve the convergence of optimization methods without loss of quality. 

Consider a single linear equation $a^\top x = b$, $x \in \mathbb{R}^n$,  where $a$ is a dense vector. The system can be equally stated as 
\begin{gather}\label{eq:conditions}
    a_i^\top x = b z^i, \quad 1\le i \le \lceil n/d\rceil, \quad \sum_{i=1}^{\lceil n/d\rceil} a_i = a, \quad \sum_{i=1}^{\lceil n/d\rceil} z^i = 1,\; z\in \mathbb{R}^{\lceil n/d\rceil}
\end{gather}
such that each $a_i$ has no more than $d$ non-zero elements. In order to guarantee that $\|a^\top x - b\|_2^2 \le \varepsilon^2$, one needs to find a vector $\psi = (x, \,z)^{\top}$, $\psi \in \mathbb{R}^{n+\lceil n/d\rceil}$, such that 
\[\sum_{i=1}^{\lceil n/d \rceil} \lceil n/d \rceil \|a_i^\top x - b z^i\|_2^2 \le \varepsilon^2\]
as under Conditions~\eqref{eq:conditions} and the Cauchy inequality one has: 
\[
    \|a^\top x - b\|_2^2 = \biggl\|\sum_{i=1}^{\lceil n/d\rceil} a_i^\top x - b z_i\biggr\|_2^2  \ge \lceil n/d\rceil^{-1} \sum_{i=1}^{\lceil n/d\rceil}\|a_i^\top x - b z_i\|_2^2. 
\]

Eqs.~\eqref{eq:conditions} correspond to a system of linear equations $\tilde A \psi = \tilde A (x, \,z)^{\top} = 0$ of the size~$(n+ d\cdot \lceil n/d\rceil) \times (n + \lceil n/d\rceil)$ with no more than $d+1$ non-zero elements in each row. 
For each $i$, row $\tilde A_i$ of $\tilde A$ is as such $\tilde A_{ij} = a_{ij}$ for $1\le j \le n$;  $\tilde A_{ij} = - b$ for $j = n+i$, and $\tilde A_{ij} = 0$ otherwise. For our convenience, we refer ${\tilde A}^x$ as the matrix consists of the first $n$ columns of $\tilde A$, and ${\tilde A}^z$ and the matrix consists of the last $n/d$ columns. 
In order to solve the problem, we apply the \texttt{NL1} algorithm, although it requires a minor correction:
\begin{gather*}
	f_\gamma(x,\,z) = \frac{1}{2}\|\tilde A 
	\psi\|_2^2 +  \frac{\gamma}{2}\sum_{i=1}^n (- x^i)_+^2 + \frac{\gamma}{2}\sum_{j=1}^{\lceil n/d\rceil} (- z^i)_+^2\to \min_{\substack{x^\top e = 1, \, x\ge 0\\ z^\top e = 1}}.
\end{gather*}

A step of the \texttt{NL1} algorithm is 
\begin{gather}\label{eq:separation}
	{\dbinom {x_{k+1}}{z_{k+1}}} = {\dbinom {x_{k}}{z_{k}}} + \argmin_{\substack{h_x:\, h^\top_x e = 0 \\ h_z: \, h^\top_z e = 0}} \left\{f_\gamma(x_k, z_k) + \nabla f_\gamma(x_k,z_k)\cdot {\dbinom {x_{k}}{z_{k}}} + \frac{L}{2} \left\|{\dbinom {h_{x}}{h_{z}}}\right\|_1^2\right\}.
\end{gather}

Notice, that Eq.~\eqref{eq:separation} is separable in variables $x$ and $z$, thus
\begin{align*}
    x_{k+1} & = x_k + \argmin\limits_{h_x: h_x^\top e = 0} \left\{\frac{1}{2} \|{\tilde A}^x h_x\| + \frac{\gamma}{2}\sum_{i=1}^n (- x^i)_+^2 + \gamma {\tilde A}^x + \gamma \sum_{i=1}^n (-x_i)_+ +  \frac{L}{2} \left\|h_{x}\right\|_1^2\right\}\\
    z_{k+1} & = z_k + \argmin\limits_{h_z: h_z^\top e = 0} \left\{\frac{1}{2} \|{\tilde A}^z h_z\| + \frac{\gamma}{2}\sum_{i=1}^{\lceil n/d\rceil} (- z^i)_+^2 + \gamma {\tilde A}^z + \gamma \sum_{i=1}^n (-z_i)_+ +  \frac{L}{2} \left\|h_{z}\right\|_1^2\right\}
\end{align*}

Recall that the optimal 
$h = (h_x,\, h_z)^{\top}$ is such as it has only two non-zero coordinates in each of $h_x$ and $h_z$, as was true for the \texttt{NL1} algorithm over the simplex, corresponding to 
\[i^+_x = \argmax_{1\le i \le n} \frac{\partial f(x_k,z_k)}{\partial x^i}\;\text{ and }\; i^-_x = \argmin_{1\le i \le n} \frac{\partial f(x_k,z_k)}{\partial x^i},\]
if 
\begin{equation*}
	\max_{1\le i \le n} \!\!\frac{\partial f_\gamma(x_k,z_k)}{\partial x^i} - \min_{1\le i \le n} \!\!\frac{\partial f_\gamma(x_k,z_k)}{\partial x^i} > \!\!\!\! \max_{1\le j \le \lceil \frac{n}{d} \rceil} \frac{\partial f_\gamma(x_k,z_k)}{\partial z^j} - \!\!\min_{1\le j \le \lceil \frac{n}{d} \rceil} \!\!\!\!\frac{\partial f_\gamma(x_k,z_k)}{\partial z^j},
\end{equation*}
while 
\[i^+_z = \argmax_{1\le j \le \lceil \frac{n}{d} \rceil} \frac{\partial f_\gamma(x_k,z_k)}{\partial z^j}\;\text{ and }\; i^-_z = \argmin_{1\le j \le \lceil \frac{n}{d} \rceil} \frac{\partial f_\gamma(x_k,z_k)}{\partial z^j},\]
otherwise. Now, for efficient implementation of the \texttt{NL1} algorithm, one needs to store the values $\partial f_\gamma/\partial z^j$ and $\partial f_\gamma/\partial x^i$ in separate binary heaps. A similar strategy can be used if sparsification of multiple rows is required. Finally, we remind the reader that single-row sparsification increases the problem dimension by $\lceil n/d\rceil$, at most.

Column sparsification is slightly more involved. Consider a function $f(x)$ in more details:
\begin{align*}
f_\gamma(x)  = & \frac{1}{2} \|Ax\|_2^2 + \frac{\gamma}{2}\sum_{i=1}^n (-x)_+^2, \qquad x = (x_1, x_2, \dots x_n)\\
{\bar f}_\gamma(x_1) = & \min_{x_{2:n}}\left\{f(x) + I_{\{x\in \Delta^n_1\}}\right\}.
\end{align*}
The function ${\bar f}_\gamma(x_1)$ is convex in $x_1$ for any fixed $\gamma > 0$ \cite[Section 3.2.5]{boyd2004convex}. That is why the value of ${\bar f}_\gamma(x_1)$ can be computed using the \texttt{NL1} algorithm while the value 
\[
f^*_\gamma = \min_{0 \le x_1\le 1}f_\gamma(x_1)
\]
requires a one dimensional binary search so that the overall time complexity $T$ of the \texttt{NL1} algorithm used to solve the PageRank problem with a single dense column~is in $O(\log (n/\varepsilon))$ higher than the one established in Theorem~\ref{thm:NL1}.

The same approach can be applied to improving the efficiency of the \texttt{S-FW} and \texttt{GK} algorithms proposed later in Sections \ref{sec:s-fw} and \ref{sec:GK} on the PageRank instances with a few dense rows or columns. 
\section{Frank--Wolfe algorithm with sparse updates}\label{sec:s-fw}

The PageRank problem, according to Eq.~\eqref{eq:pgm-setup}, can be stated as 
\begin{gather*}
f(x) = \frac{1}{2} \|Ax\|_2^2 \to \min_{x \in \Delta^n_1}.
\end{gather*}
We use the Frank--Wolfe conditional gradient \cite{frank1956algorithm,jaggi2013revisiting} to solve the problem above.

We choose the starting point $x_0$ of the algorithm in {an arbitrary vertex} of the unit simplex. Then on each step we solve 
\begin{gather}\label{eq:fw-lin}
h_k^\top \nabla f(x_k) \to \min_{y \in \Delta^n_1}.
\end{gather}

Furthermore the solution $y_k$ of Eq.~\eqref{eq:fw-lin} has  
only one non-zero coordinate~$y_k^{i_k}$, corresponding to 
\[i_k = \argmin_{1\le i\le k} \partial f(x)/ \partial x^i.\]

Then, the update rule is
\begin{gather*}
x_{k+1} = (1-\gamma_k) x_k + \gamma_k h_k, \quad \gamma_k = \frac{2}{k+1}, \quad k\ge 1.
\end{gather*}

According to \cite{nesterov2013introductory}, $f(x_k) - f(x)$ is bounded from above as:
\begin{gather}\label{eq:fw-it}
f(x_k) - f^* = f(x_k) \le \frac{2 L_1 \max_{x,y \in \Delta^n_1}\|x - y\|_1^2}{k+1} \le \frac{8 L_1}{k+1},
\end{gather}
where $L_1^2 = \max_{x \in \Delta_1^n} \|Ax\|_2^2 \le 2.$ Thus, in order to guarantee $f(x_k) \le \varepsilon^2/2$, one needs at most $32\varepsilon^{-2}$ iterations. 

Consider a step of the algorithm in more detail. Denote $\beta_k$ as 
\[\beta_k = \prod_{r=1}^{k-1} (1-\gamma_r), \; z_k = x_k/\beta_k, \; \tilde{\gamma}_k = \gamma_k/\beta_{k+1}, {\text{with } \beta_0 = 1}.\]
A step of the Frank--Wolfe algorithm is
\[z_{k+1} = z_k + \tilde{\gamma}_k h_k, \; {\text{with }z_1 = x_1}.\]
where $h_k$ is a solution of Eq.~\eqref{eq:fw-lin}.
 Moreover, only one coordinate of~$h_k$: 
\[
{i_k} = \argmin_{1\le i\le k} \partial f(x_k)/\partial x^{i} = \argmin_{1\le i\le k} A^\top A z_k^i
\]
is other than zero. Therefore, the update of the \texttt{S-FW} algorithm is similar to the Gauss-Southwell rule studied in detail in \cite{nutini2015coordinate} for minimization of the strongly convex functions. 
\begin{algorithm}[!t]
\DontPrintSemicolon 
\KwIn{$d$-sparse transition matrix $P$, starting point $x_0$ in one of the vertices of unit simplex and $A = P^\top - I$; Set $k=1$ and $\beta_0 = 1$}
\KwOut{$x_k:\; {\|Ax_k\|_2} \le \varepsilon$}
  \While{$f(z_k) > \beta_k^2 \varepsilon^2$} {
     $i_k = \argmin_{1\le i\le n} {\partial f(x_k)}/{\partial x^i}$\;
     {Set $\gamma_k = 2/(k+1)$, $\beta_{k+1} = \beta_k (1-\gamma_{k})$, and $\tilde \gamma_k = \gamma_k/\beta_{k+1}$}\;
     Update step direction $h_k = \tilde \gamma_k \cdot \delta$, where $\delta^{i_k} = 1$, and $\delta^i = 0$, $i \neq i_k$\;
     Update argument: $z_{k+1} = z_k + \tilde \gamma_k h_k$\;
     Update gradient: $\frac{\nabla f(x_{k+1})}{\beta_k} = \frac{\nabla f(x_{k})}{\beta_k} + \tilde \gamma_k A^\top A h_k$\;
     Update function value $f(x_{k+1})$: \, $\frac{f(x_{k+1})}{\beta_k^2} = \frac{f(x_{k})}{\beta_k^2} + \tilde \gamma_k h_k^\top \frac{\nabla f(x)}{\beta_k} + {\|Ah_k\|_2^2\tilde\gamma_k^2}/2$\;
     $k = k+1$\;
  }
\Return{$x_k = z_k \beta_k$}\text{\quad //Compute $x_k$ on the last iteration only}\;
\caption{{\sc \texttt{S-FW}: Frank--Wolfe algorithm for PageRank}}
\label{algo:fw}
\end{algorithm}

Using the doubly-linked list of binary heaps, described in Section~\ref{sec:NL1}, the minimal coordinate of $A^\top A z_k$ can be computed in $O(d^2 \log (2 + n/d^2))$ time. Then $A^\top A z_{k+1}$ is 
\begin{gather}\label{eq:upd-fw}
	A^\top A z_{k+1} = A^\top A z_k + \tilde{\gamma}_k A^\top A y_k,
\end{gather}
and 
\begin{gather*}
\frac{\nabla f(x_{k+1})}{\beta_k} = \frac{\nabla f(x_{k})}{\beta_k} + \frac{\tilde \gamma_k \beta_k A^\top A h_k}{\beta_k},
\end{gather*}
and also
\begin{gather*}
\frac{f(x_{k+1})}{\beta_k^2} = \frac{f(x_{k})}{\beta_k^2} + \tilde \gamma_k h_k^\top \frac{\nabla f(x)}{\beta_k} + {\|Ah_k\|_2^2\tilde\gamma_k^2/2}.
\end{gather*}

For a $d$-sparse matrix $A$, the time required to compute $A^\top A z_{k+1}$ from $A^\top A z_{k}$ is $O(d^2\log(2 + n/d^2))$, as well. One can compute $x_k$ having $z_k$ as $x_k = \beta_k x_k$ in $O(n)$ time.

Combining Eqs.~\eqref{eq:fw-it} and \eqref{eq:upd-fw}, we have the following complexity estimate for the algorithm.  

\begin{theorem}\label{thm:S-FW}
Algorithm \ref{algo:fw} requires at most $k = 32 \varepsilon^{-2}$ iterations to guarantee $\|P^\top x - x\|_2 \le \varepsilon$ for any $d$-sparse transition matrix $P$. The overall time complexity of the algorithm does not exceed 
\[
T = O\left(n + \frac{d^2 \log (2 + n/d^2)}{\varepsilon^2}\right). 
\]
\end{theorem}

\textit{Discussion.} Theorem~\ref{thm:S-FW} implies sub-linear convergence in the number of non-zero elements of the transition matrix $P$. This is not surprising, since we assume that the underlying graph structure, along with required smoothness, are known a priori. 

It remains an open question for the authors to improve convergence rate in terms of the accuracy and maximal degree of the transition graph while preserving linear dependence in the problem dimension.

\section{Saddle point setup for PageRank}\label{sec:GK}

In the PageRank problem, one often requires an accurate approximation of a few of the largest coordinates representing the most relevant websites rather than the full PageRank vector. To this end, we propose an algorithm to approximate the PageRank vector in $\ell_\infty$-norm. 
Below, we consider the problem
\begin{gather*}
f(x) = \|(P^\top - I)x\|_\infty  = \|Ax\|_{\infty }\to \min\limits_{x \in \Delta_1^n}.
\end{gather*}

Following \cite{nemirovski2009robust}, we set up the problem as
\begin{gather}\label{eq:matrix-game}
\min_{x\in \Delta^n_1} \max_{\|{\tilde{y}}\|_1 \le 1}  \langle Ax, {\tilde{y}} \rangle 
       =  \min_{x\in \Delta_1^n}\max_{y \in \Delta^{2n}_1}
        \langle Ax, Jy \rangle 
        = \min_{x\in \Delta^n_1}\max_{y \in \Delta^{2n}_1} \langle x, \tilde A y \rangle,
\end{gather}
where $\tilde A = A^\top J$, $J  = [I_n, -I_n]$, and $I_n$ is the $n \times n$ identity matrix. 

We propose a sub-linear-time algorithm to approximate a bilinear matrix game representing the PageRank, Problem~\eqref{eq:matrix-game}. Let $\tilde A_{ij}$ be a gain for Player $A$ (loss of player $B$), if Player $A$ plays strategy $i$ and $B$ plays strategy $j$, $1 \le i \le n$, and $1\le j \le 2n$. Consider the loss function for Player $B$ at step $k$:
\[
f(x, y_k) = x^\top \tilde A y_k, \quad x \in \Delta_1^n, 
\]
where $y_k \in \Delta_1^{2n}$ is a vector with a single non-zero coordinate corresponding to the strategy of Player $A$.  
We also emphasize that $y_k$ depends 
on the whole history of the game. 
Let $C$ be the cost of the matrix game:
\begin{equation}
    \label{eq_saddle_opt}
C = \max_{y \in \Delta^{2n}_1} \min_{x \in \Delta^n_1} y^\top \tilde A x = \min_{x \in \Delta^n_1} \max_{y \in \Delta^{2n}_1} y^\top \tilde A x = \min_{x\in \Delta_1^n} \|Ax\|_\infty = 0, 
\end{equation}
and 
\begin{gather}\label{eq:gk-bounds}
\min_{x \in \Delta_1^n} \frac{1}{N} \sum_{i=1}^N f(x, y_i) \ge C\ge \max_{y \in \Delta_1^{2n}}  \frac{1}{N} \sum_{i=1}^N f(x_i, y),
\end{gather}
for any sequences $\{x_i\}_{i=1}^N$, $\{y_i\}_{i=1}^N$ if for any $i$: {$x_i \in \Delta_1^n$, $y_i \in \Delta_1^{2n}$}, $1 \le i \le n$. In the subsequent of the section, we consider $\{(x_i, y_i)\}_{i=1}^N$ with a single non-zero coordinate each. 

\begin{algorithm}[!t]
\DontPrintSemicolon 
\KwIn{(unnormalized) probability distribution $p$ given by a binary tree with leafs values $p_i$ such that the 
probability of any leaf $j$ is $p_j/\sum_{i=1}^n p_i$,  $p_i > 0$ for any $i:\; 1\le i \le n$, $\sum_{i=1}^n p_i > 0$, with and update rule of coordinate $k$ according to~Eq.~\eqref{eq:prob-upd} 
\begin{gather*}
	{\tilde p}_k \propto {p}_k \exp\left(- \psi_k\right),
\end{gather*}
}
\KwOut{$x\sim \tilde p$, a sample $x$ follows the updated distribution $\tilde p$}
// Distribution update\;
$u = k$ // Start with the leaf $k$\;
  \While{$u \neq \text{root}$} {
     $p_u \leftarrow p_u + p_k (\exp\left(- \psi_k\right)-1)$;\;
     $u \leftarrow $ parent of $u$\;
  }
// Sampling $x\sim \tilde p$\;
$u = k$ // Start with the leaf $k$\;
  \While{$u \neq \text{leaf}$} {
   	 Let $\nu, \omega$ be children of $u$
     \[u = \begin{cases}
     \nu, & \text{ with probability } p_{\nu}/(p_{\nu} + p_{\omega}), \\
     \omega, & \text{ otherwise}.
     \end{cases}\]
  }
\Return{$u$}\;
\caption{{\sc update rule for a probability distribution}}
\label{algo:upd}
\end{algorithm}

\begin{algorithm}[!t]
\DontPrintSemicolon 
\KwIn{$d$-sparse transition matrix $P$, starting point $x_0$ in one of the vertices of the unit simplex, and learning rate $\gamma$}
\KwOut{$x_k:\; f({\bar x}_k, {\bar y}_k) = {\bar x}_k \tilde A {\bar y}_k \le {\varepsilon}$,  implies $\|P^\top x_k - x_k\|_\infty \le \varepsilon$}
$\pi = ({(2n)}^{-1}, \dots, {(2n)}^{-1})$, \quad $p = (n^{-1}, \dots, n^{-1})$, \; starting point $(x_0, y_0)$ in one of the vertices of $\Delta_1^n \times \Delta_1^{2n}$\;
  \While{$f({\bar x}_k, {\bar y}_k) > {\varepsilon}$} {
  // Player $A$ turn\;
     Choose at random $i_k^A$, such that $\mathbb{P}(i_k = j) = \pi_i$;\;
     Assume $y^{i_k^A}_k = 1$, and $y_k^i = 0$ if $i\neq i_k^A$;\;
     Update $\pi_k$// see Algorithm~\ref{algo:upd} for details
     \[
     \pi^i_{k+1} \propto \pi^i_k \exp 
     	\left(
             {\gamma_y}
            {\tilde A}_{i, j^B_k}
        \right)
     \]
     // Player $B$ turn\;
     Choose at random $j_k^B$, such that $\mathbb{P}(j_k^B = j) = p_j$;\;
     Assume $x^{j_k^B}_k = 1$, and $x_k^j = 0$ if $j\neq j_k^B$;\;
     Update $p_k$ // see Algorithm~\ref{algo:upd} for details
     \[
     p^j_{k+1} \propto p^j_k \exp 
     	\left(
            - {\gamma_x} 
            {\tilde A}_{i_k^A, j}
        \right)
     \]\;
     \vskip -4mm
     {Update an average point $({\bar x}_k, {\bar y}_k)$. Indeed, more time-efficient is to update $k \cdot({\bar x}_k, {\bar y}_k)$ as this involves only sparse operations according to:}
     \[k \cdot {\bar x}_k = \sum\limits_{t = 1}^k x_t = (k-1){\bar x}_{k-1} + x_k, \quad {k \cdot} {\bar y}_k = \sum\limits_{t = 1}^{k} y_t = (k-1){\bar y}_{k-1} + y_k\]
     {an can be done in $O(d)$ time}\;
     {Update the function value, $f({\bar x}_k, {\bar y}_k)$, using sparse operations only:}
     \begin{align*}
     k^2 f({\bar x}_k, {\bar y}_k) & = (k{\bar x}_k)^\top \tilde A (k{\bar y}_k) \\
     & = 
     (k-1)^2 f({\bar x}_{k-1}, {\bar y}_{k-1})
     + (k-1)\bar x_{k-1}\tilde A y_k + x_k \tilde A ((k-1){\bar y}_k) + x_k \tilde A y_k
     \end{align*}
     as $k{\bar x}_k = (k-1) {\bar x}_{k-1} + x_k$, and $k{\bar y}_k = (k-1) {\bar y}_{k-1} + y_k$. {The update requires $O({d})$ time. }
  }
\Return{$({\bar x}_k, {\bar y}_k)$}\;
\caption{{\sc $\ell_\infty$~approximation to the PageRank problem}}
\label{algo:GK}
\end{algorithm}

To solve the problem, we assume the following randomized strategy for the Player $B$ played against any strategy of the Player $A$:
\begin{enumerate}
\item Let $p_1 = (n^{-1}, \dots, n^{-1})$; 
\item Choose at random $j_k$, such that $\mathbb{P}(j_k = j) = p_{j}^k$;
\item Assume $x_{j_k}^k = 1$ and $x_j^k = 0$ for all $j\neq j_k$;
\item Update 
\begin{gather}\label{eq:prob-upd}
	p_j^{k+1} \propto p_j^k \exp\left(- {\gamma_x} \tilde A_{i_k, j}\right),
\end{gather}
\noindent where $i_k$ is a strategy that Player $A$ chooses at step $k$. 
\end{enumerate}

The crucial gain in the efficiency of the algorithm is due to time-efficient updates at stage 3. Indeed, consider a binary tree with its leaves corresponding to the variables, and constructed in such a way that the value $p_v$ assigned to a node $v$ is a total probability of all leaves having $v$ as a predecessor. If we update the weight of a leaf according to Eq.~\eqref{eq:prob-upd}, we also update each vertex $u$ belonging to the path from the leaf to the root as 
\[
	p_u = p_u + \xi, \quad \xi = p_j^{k+1} - p_j^k.
\]

In order to sample $x \sim p^{k+1}$, we start from the root of the tree and proceed to its child~$a$ with probability $p_a/(p_a + p_b)$. Otherwise, we proceed to its sibling $b$, where $p_a$, and $p_b$ are the values assigned to $a$ and $b$, respectively. We repeat the same procedure for each node one a path from the root to one of the leafs of the tree.  Algorithm~\ref{algo:upd} formalizes this argument. 

Using the same strategy for the Player $A$, we establish the convergence rate to the Nash equilibrium $(x^*,\, \omega^*)$, which  solves Problem~\eqref{eq:matrix-game} in Theorem~\ref{thm:GK-ch-p}. Algorithm~\ref{algo:GK} contains all necessary details. It is worth mentioning an interpretation of the algorithm, e.g., on each iteration it make an update following to a sparse projection  of the gradient to the unit simplex according to the KL divergence.

\begin{theorem}\label{thm:GK-ch-p}
Algorithm \ref{algo:GK} after {$N\ge 4\varepsilon^{-2}\Big(\ln(2n)+ \ln(n) + 16\ln(1/\delta) \Big)$ iterations with a constant step-size  $\gamma_x = \sqrt{2(\log n)/N}$ and $\gamma_y = \sqrt{2(\log (2n))/N}$}  results in a point 
$
(\bar x_N, \bar y_N)$ 
such that, with probability at least $1 - \delta$, for any $\delta > 0$ one has:
\[
    {\|A \bar x_N\|_\infty \le  \varepsilon.}
\]
Moreover, the total running time of the algorithm is bounded from above as 
\[T = O\left(n + \frac{d \log n \log \frac{n}{\delta}}{\varepsilon^2}\right).\]
\end{theorem}

\noindent The proof of the theorem is provided in Appendix \ref{sec:appB}. 

\textit{Discussion. The mirror descent and the dual averaging perspectives.}  The proposed algorithm is essentially a mirror descent with randomized projection of the gradient on a unit simplex according to KL divergence.
Let us consider the problem of minimizing of the left hand-side of Eq.~\eqref{eq:gk-bounds} in more details: 
\begin{gather}\label{eq:soo}
\frac{1}{N}\sum_{i=1}^N f(x, y_i) = \frac{1}{N}\sum_{i=1}^N 
 \langle x, A^\top J y_i \rangle  \to \min_{x\in \Delta_1^n},
\end{gather}
where $\{y_i\}_{i=1}^n$ is a sequence of unit coordinate vectors. Denote $f_i(x) \doteq f(x, y_i)$ for $i \ge 1$. 
Recall the setup of the mirror descent algorithm \cite{nemirovski2009robust}. Let 
$
\omega(x) = \sum_{i=1}^n x^i \log x^i 
$
be the distance-generating function, which is is 1-strongly convex with respect to the $\ell_1$ norm. A step of the dual averaging algorithm \cite{nesterov2009primal} with step-size ${\gamma_x}$ is:
\begin{gather}\label{eq:rand-upd}
\quad {z}_{k} = {z}_{k-1} - {\gamma_x} \nabla f_k(x_{k}), \quad x_{k+1} = \nabla \omega^* ({z}_{k}), 
\end{gather}
where
$\omega^*({z}) = \sup_{x\in \Delta^n_1} \left\{{z}^\top x - \omega(x)\right\} = \log \left\{\sum_{j=1}^n \exp({z}_j)\right\}$. An update of $x_{k+1}$  in Eq.~\eqref{eq:rand-upd} can be also viewed as a projection of ${z}_{k+1}$ to the unit simplex in accordance with Kullback-Leibler divergence. Indeed  Eq.~\eqref{eq:rand-upd} is a step of mirror descent algorithm for simplex constrained problems as well \cite[Appendix A]{allen2014linear}, \cite{banerjee2005clustering}, \cite{juditsky2011first}. 
A randomized version of the update is then 
\begin{gather}\label{eq:upd-mp}
x_{k+1} = e_i, \text{ with probability } p_{k+1}^i = \frac{p^i_k \exp\left( - {\gamma_x} \frac{\partial f_k (x_k)}{\partial x^i} \right)}{\sum_{t=1}^n p_k^t\exp\left( -  {\gamma_x} \frac{\partial f_k (x_k)}{\partial x^t} \right)},
\end{gather}
where $e_i$ is a unit vector with a single non-zero coordinate corresponding to index $i$. Since $f(x,y) = y^\top \tilde A x$, update \eqref{eq:upd-mp} is the same as the update in~Algorithm~\ref{algo:GK}. 

\section{Implementation details and case study}\label{sec:case}

All algorithms proposed in the paper are implemented in \texttt{C++}. We test our code with different versions of \texttt{GCC} (GNU Compiler Collection), \texttt{clang} (C language family front-end for LLVM), and \texttt{icc} (Intel C
Compiler) compilers under GNU/Linux, Microsoft Windows, and Mac OS X. 
We conduct the experiments using:
\begin{itemize}
 \item Ubuntu server 16.04.6 LTS, x86\_64
 \item Intel Core i5-2500K, 16 Gb RAM
 \item \texttt{GCC}-5.4.0 to compile \texttt{C++} code,
\item Assembly parameters: \texttt{-std=c++11 -O2 -mcmodel=small -DNDEBUG}
\end{itemize}

We test our algorithms in different dimensions using the following three test beds:
\begin{enumerate}
\item $d$-diagonal matrix for $n_d = 1,3,5, \dots$. Each row/column of these matrices contains $(n_d - 1)/2 + 1 \le d \le n_d$ non-zero elements; 
\item randomly-generated matrices with $d$ non-zero elements (on average);
\item and web-graphs from the Stanford University graph collection \footnote{\url{http://snap.stanford.edu/data/\#web}}
\end{enumerate}

We use accuracy $\varepsilon = 10^{-4}$ in each of our experiments; $x_0 = (1, 0, \dots, 0)$ is used as a starting point for the \texttt{NL1} and \texttt{S-FW} algorithms, and we terminate the algorithms if $f(x_k) = \|Ax\|_2^2/2 \le \varepsilon^2/2$. Computational time reported for the case study includes time required by optimization method, $A$ and $A^\top$ generating times, an initialization of all data structures used, as well as initial gradient/function computation. 

\begin{table}[h!]
 \begin{center}
  \begin{tabular}{lr||r|r|r|r|r}
  \multirow{3}{*}{web-graph} &
  \multirow{3}{*}{} &
    \multicolumn{5}{c}{\# non zeros} \\ \cline{3-7}
    & & \multicolumn{2}{c|}{in a row} & \multicolumn{2}{c|}{in a column} & \multirow{2}{*}{average} \\ \cline{3-6}
    & & min & max & min & max\\ 
  \hline

  Stanford,  & $n$ = 281 903 & 2 & 38 607 & 1 &   256 &  9.20  \\
  NotreDame, & $n$ = 325 729 & 2 & 10 722 & 1 & 3 445 &  5.51 \\
  BerkStan,  & $n$ = 685 230 & 1 & 84 209 & 1 &   250 & 12.09 \\
  Google,    & $n$ = 875 713 & 1 &  6 327 & 1 &   457 &  6.83 \\
  \end{tabular}
  \end{center}
  \caption{Structure of the matrix $A$ for the graphs from the Stanford web-graph collection. Columns of the transition matrix are more dense than the rows and on average both columns and rows contain a few non-zero elements only.}
  \label{tbl:a_webgraphs_char}
\end{table}

The numerical experiments described below allow the following conclusions to be drawn:

\begin{enumerate}
\item 
In our experiments (see Figure~\ref{fig:gk}) algorithm \texttt{GK} converges  sufficiently fast for the desired precision. But after a large number of iterations, the value of $f(\bar x)$ starts to grow, and the resulting point does not satisfy the accuracy condition. The reason for this is that the values of several probabilities $p_i$ become extremely large and out of range after a number of iterations. That leads to significant numerical errors in estimating residual small probabilities and overall unsatisfactory performance of the algorithm. Rescaling the probability vectors does not change the behavior of the algorithm. The described effect decreases for larger $n$; refer to Fig.~\ref{fig:gk} for details.

 \begin{figure}[h!]
  \begin{center}
    \includegraphics[width=0.95\linewidth]{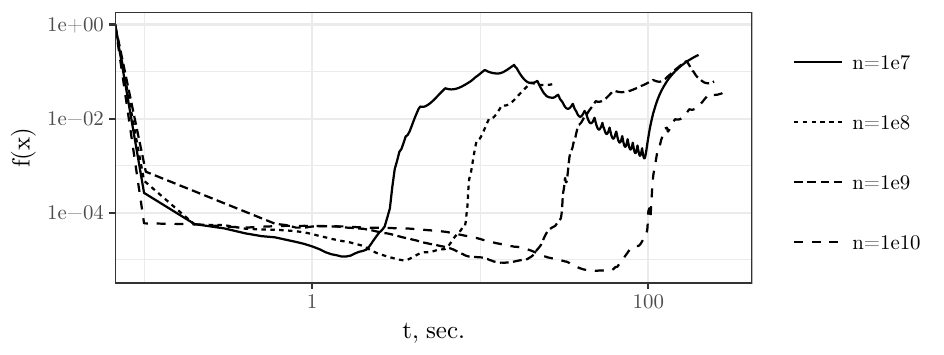}   
  \end{center}
  \caption{Convergence of the \texttt{GK} algorithm for various dimensions $n$. $A$ is a random $n\times n$ matrix, with a number of non-zeros in each row and column $d=3$. Practical performance of the \texttt{GK} algorithm is limited due to unavoidable errors in estimating small probabilities. }
  \label{fig:gk}
 \end{figure}
 
To summarize, the theoretical bounds for the \texttt{GK}   algorithm differ markedly from from those seen in practical performance. 
  
\item Computational time for $d$-diagonal matrices $A$ is much smaller than that for random matrices 
(see Figures~\ref{fig:a_diag} and \ref{fig:a_rand} for details).
This is due to the fast cache operations, which require far fewer memory reads for sequential data. Updates to computational trees/heaps are also performed in sequential elements, which improves time performance as well. Also, this permits the dramatic reduction of the dependence on the actual problem dimension in practice.  
Thus in~Table~\ref{tbl:a_diag} for $n_d=3$ and accuracy $\varepsilon = 10^{-4}$ the computational time has increased less than twice for $n=10^{8}$ compared with $n = 10^2$.

Conversely, for the random matrices, caching does not give the same improvement in speed. This significantly decreases the actual performance of the algorithms; see Table~\ref{tbl:a_random} for details. 

\begin{table}[h!]
  \begin{center}
  \begin{tabular}{r||r|r||r|r}
  & \multicolumn{2}{c||}{\texttt{NL1}} & 
    \multicolumn{2}{c  }{\texttt{S-FW}} \\ \cline{2-5}
  $n$ & time, sec. & iteration & time & iterations \\
  
  \hline
  \multicolumn{3}{l}{$n_d = 3$; $2 \le d \le 3$} \\
  \hline
  
  $10^2$ & 4.089 & 3 948 632 & 0.007 & 14 142 \\
  $10^3$ & 4.221 & 3 950 392 & 0.008 & 14 142 \\
  $10^4$ & 4.575 & 3 950 392 & 0.009 & 14 142 \\
  $10^5$ & 4.814 & 3 950 392 & 0.010 & 14 142 \\
  $10^6$ & 5.143 & 3 950 392 & 0.010 & 14 142 \\
  $10^7$ & 5.566 & 3 950 392 & 0.010 & 14 142 \\
  $10^8$ & 6.021 & 3 950 392 & 0.010 & 14 142 \\ 

  \hline  
  \multicolumn{3}{l}{$n_d = 11$; $6 \le d  \le 11$} \\
  \hline
  
  $10^2$ & 14.655 & 2 100 964 & 0.041 & 14 749 \\
  $10^3$ & 37.796 & 5 101 072 & 0.041 & 16 956 \\
  $10^4$ & 39.170 & 5 101 072 & 0.062 & 19 995 \\
  $10^5$ & 39.897 & 5 101 072 & 0.064 & 24 495 \\
  $10^6$ & 41.004 & 5 101 072 & 0.065 & 24 495 \\
  $10^7$ & 43.917 & 5 101 072 & 0.068 & 24 495 \\ 
  
  \hline  
  \multicolumn{3}{l}{$n_d = 51$; $26 \le d \le 51$} \\
  \hline
  
  $10^3$ & 529.240 & 5 216 119 & 1.552 & 46 447 \\
  $10^4$ & 535.348 & 5 216 119 & 1.045 & 29 991 \\
  $10^5$ & 537.419 & 5 216 119 & 1.741 & 49 235 \\
  $10^6$ & 549.782 & 5 216 119 & 1.758 & 49 235 \\
  $10^7$ & 552.271 & 5 216 119 & 1.789 & 49 235 \\
  
  \hline
  \multicolumn{3}{l}{$n_d = 101$; $51 \le d \le 101$} \\
  \hline
  
  $10^4$ & 1 935.198 & 5 175 085 & 6.464 & 49 925 \\
  $10^5$ & 1 962.307 & 5 175 085 & 9.097 & 68 646 \\
  $10^6$ & 1 940.331 & 5 175 085 & 9.134 & 68 646 \\
  
  \end{tabular}
  \end{center}

  \caption{Time in seconds required to solve the PageRank problem. $A$ is a $d$-diagonal matrix. The \texttt{S-FW} algorithm outperforms the \texttt{NL1} algorithm for most of the instances and have better scalability with the dimension of the problem. }
 \label{tbl:a_diag}
\end{table}

\begin{figure}[h!]
  \begin{center}
  \includegraphics[width=1.0\linewidth]{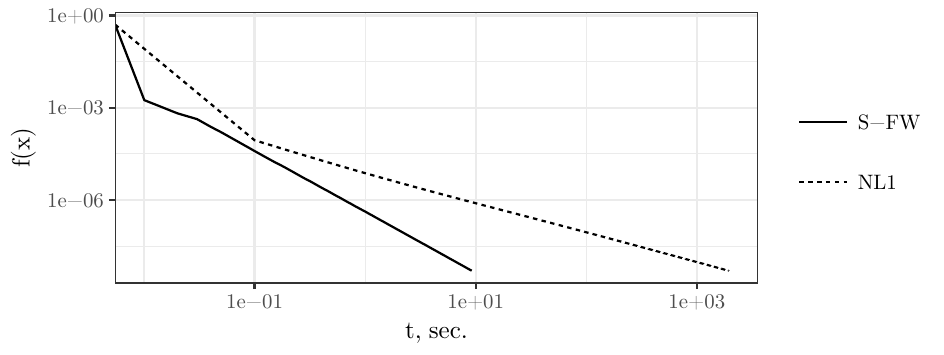}
  \end{center}
  \caption{Computational time for the PageRank problem, $A$ is a $n_d$-diagonal matrix, dimension $n=10^6$, number of non-zero diagonals $n_d=101.$ The \texttt{S-FW} algorithm significantly outperforms the \texttt{NL1} algorithms for large scale problems.}
  \label{fig:a_diag}
\end{figure}

\begin{table}[h!]
 \begin{center}
  \begin{tabular}{r||r|r||r|r}
  & \multicolumn{2}{c||}{\texttt{NL1}} & 
    \multicolumn{2}{c  }{\texttt{S-FW}} \\ \cline{2-5}
  $n$ & time & iterations & time & iterations \\
  
  \hline
  \multicolumn{3}{l}{$d = 3$} \\
  \hline
  
  $10^2$ &   0.003 &      1 999 &   0.023 &     39 734 \\
  $10^3$ &   0.031 &     17 748 &   0.118 &    190 601 \\
  $10^4$ &   0.233 &    141 739 &   0.414 &    632 954 \\
  $10^5$ &   2.374 &    840 617 &   2.107 &  2 009 854 \\
  $10^6$ &  16.171 &  4 020 388 &   9.355 &  6 203 826 \\
  $10^7$ &  56.694 & 11 669 495 &  32.442 & 17 916 520 \\
  $10^8$ & 173.070 & 19 988 053 & 121.258 & 43 390 838 \\
  
  \hline  
  \multicolumn{3}{l}{$d = 11$} \\
  \hline
  
  $10^2$ &   0.013 &        590 &   0.173 &     44 706 \\
  $10^3$ &   0.072 &      5 106 &   0.593 &    142 109 \\
  $10^4$ &   0.568 &     40 029 &   2.123 &    450 873 \\
  $10^5$ &   6.342 &    299 382 &  10.374 &  1 482 735 \\
  $10^6$ &  78.383 &  2 025 423 &  60.715 &  4 753 809 \\
  $10^7$ & 503.385 & 11 272 158 & 219.988 & 14 693 667 \\ 
  
  \hline  
  \multicolumn{3}{l}{$d = 51$} \\
  \hline
  
  $10^3$ &      0.891 &      3 851 &    11.681 &    162 015 \\
  $10^4$ &      8.383 &     31 372 &    42.824 &    510 444 \\
  $10^5$ &     77.137 &    241 191 &   164.751 &  1 621 686 \\
  $10^6$ &  1 300.194 &  1 683 845 & 1 152.805 &  5 082 774 \\
  $10^7$ & 11 250.461 & 10 627 974 & 5 432.107 & 17 479 622 \\
  
  \hline
  \multicolumn{3}{l}{$d = 101$} \\
  \hline
  
  $10^4$ &    29.540 &    29 127 &   168.124 &   529 685 \\
  $10^5$ &   304.419 &   225 146 &   650.878 & 1 696 708 \\
  $10^6$ & 4 692.729 & 1 607 834 & 4 619.220 & 5 267 738 \\
  
  \end{tabular}
 \end{center}
 \caption{Time in seconds required to solve the PageRank problem. $A$ is a random matrix. The \texttt{NL1} algorithm outperforms the \texttt{S-FW} algorithm in most of the sparse and low-dimensional instances, while the \texttt{S-FW} algorithm is preferable for large scale cases.}
 \label{tbl:a_random}
\end{table}

\begin{figure}[h!]
  \begin{center}
  \includegraphics[width=\linewidth]{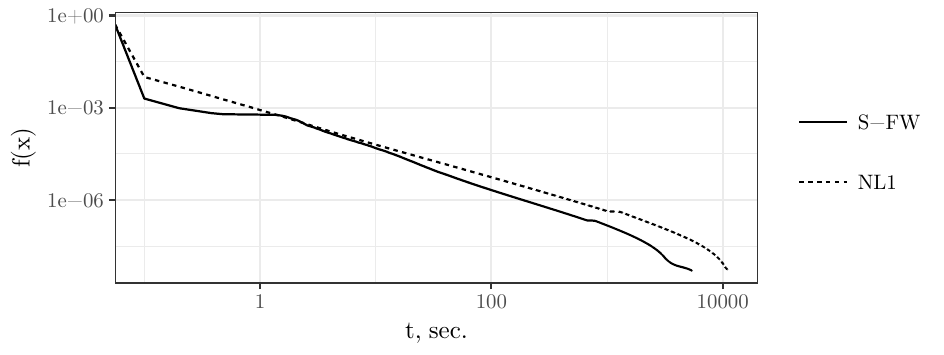}  
  \end{center}

  \caption{Time complexity for the PageRank problem. $A$ is a random matrix, dimension $n=10^7$, average number of non-zeros in each row and column $d=51$. The \texttt{NL1} and \texttt{S-FW} algorithms have almost the same computational time.}
  \label{fig:a_rand}
\end{figure}

\item Surprisingly, the time complexity of the \texttt{S-FW} algorithm for the Stanford web-graph collection is much less than that for the \texttt{NL1} algorithm (see Table~\ref{tbl:a_webgraphs}). Unfortunately, for two problems on the list, the \texttt{NL1} algorithm performance is not sufficiently high. We propose that this is due to the fact that the \texttt{NL1} algorithm modifies two variables per iteration and often involves very expensive, dense updates compared with to the \texttt{S-FW} algorithm (see~Table~\ref{tbl:a_webgraphs_char} for the information about the sparsity of the transition matrices). Table~\ref{tbl:r_sc_berkstan} contains information about the average iteration complexity of the \texttt{NL1} algorithm; it is much higher than that for the \texttt{S-FW} algorithm, which supports our conjecture,  particularly for the web-BerkStan dataset. Recall that this property is true without additional matrix sparsification (refer to Section \ref{sec:NL1}).

The performance of the \texttt{S-FW} and \texttt{NL1} algorithms is shown in 
Figures~\ref{fig:web_stanford}--\ref{fig:web_google_iter}.

\begin{table}[h!]
 \begin{center}
  \begin{tabular}{l|r||r|r||r|r}
  & & \multicolumn{2}{c||}{ \texttt{NL1}} & 
    \multicolumn{2}{c  }{ \texttt{S-FW}} \\ \cline{3-6}
  web-graph & $n$ & time, sec. & iterations & time, sec. & iteration \\
  \hline
  
  Stanford  & 281 903 &      0.145 &     93 152 & 0.008 & 14 142 \\
  NotreDame & 325 729 &    700.810 &  3 816 436 & 0.526 & 38 014 \\
  BerkStan  & 685 230 & 38 161.847 & 12 315 700 & 0.536 & 19 990 \\
  Google    & 875 713 &    113.643 &  1 083 996 & 0.278 & 37 313 \\
  
  \end{tabular}
  \end{center}
  \caption{Time in seconds required to solve the PageRank problem for web-graphs from the Stanford graph collection. The \texttt{S-FW} algorithm achieves significantly better time performance compared to the \texttt{NL1} algorithm. }
  \label{tbl:a_webgraphs}
\end{table}

\begin{table}[h!]
 \begin{center}
  \begin{tabular}{l||r|r||r|r}
  & \multicolumn{2}{c||}{Stanford, $n = 281903$} &
    \multicolumn{2}{c  }{Google, $n = 875 713$} \\ \cline{2-5}
  method & time, sec. & iterations & time, sec. & iteration \\
  \hline
  
  \texttt{S-FW}     &   0.008 & 14 142 &      0.278 &    37 313 \\
  \texttt{S-FW(NS)} &  75.438 & 14 142 &    451.131 &    38 672 \\
  \texttt{NL1}      &   0.145 & 93 152 &    113.643 & 1 083 996 \\
  \texttt{NL1(NS)}  & 458.493 & 93 152 & 13 507.423 & 1 220 868 \\
  \texttt{FGM(NS)}   &  82.978 & 12 464 &    423.008 &    22 811 \\
  \end{tabular}
  \end{center}
  \caption{Time in seconds required to solve the PageRank problem for web-graphs with sparse and non-sparse (NS) versions of \texttt{S-FW} and \texttt{NL1} methods. We compare sparse versions of the algorithms with the non-sparse ones (\texttt{NS}). We have compared our algorithms with the state-of-the-art Similar Triangles Algorithm \cite{gasnikov2018universal,nesterov2018lectures},
  \label{tbl:a_webgraphs_ns}, 
  which is essentially an extension to the Fast Gradient method (FGM).}
\end{table}

\begin{table}[h!]
 \begin{center}
  \begin{tabular}{r|l||r|r||r|r}
  \multicolumn{2}{c||}{} & \multicolumn{2}{c||}{Stanford} & \multicolumn{2}{c}{BerkStan} \\ \cline{3-6}
  \multicolumn{2}{c||}{\# elements} & NL1 & FW & NL1 & FW \\
  \hline
  
  \multirow{3}{*}{$d_r$}
  & min    &  1.0 & 1.0 &     1.0 &     1.0 \\
  & max   & 34.0 & 4.0 & 84 209.0 & 84 209.0 \\
  & average &  3.9 & 3.9 &  2 278.4 &   148.6 \\
  \hline

  \multirow{3}{*}{$d_c$}
  & min    &  2.0 & 2.0 &   1.0 &  1.0 \\
  & max  & 37.0 & 3.0 & 244.0 & 83.0 \\
  & average &  2.9 & 2.8 &  15.7 &  6.2 \\
  \hline

  \multirow{3}{*}{$d_r \cdot d_c$}
  & min    &    3.0 &  3.0 &        2.0 &       2.0 \\
  & max   & 1 258.0 & 12.0 & 15 494 456.0 & 6 989 347.0 \\
  & average &   11.7 & 11.3 &    84 304.3 &    7 507.5 \\
 \end{tabular}
 \end{center}
 \caption{Iteration complexity for the Stanford graph collection. The \texttt{S-FW} algorithm has performed a much fewer number of updates, resulting in a higher performance compared to the \texttt{NL1} algorithm.}
 \label{tbl:r_sc_berkstan}
\end{table}

\begin{figure}[h!]
  \begin{center}
  \includegraphics[width=\linewidth]{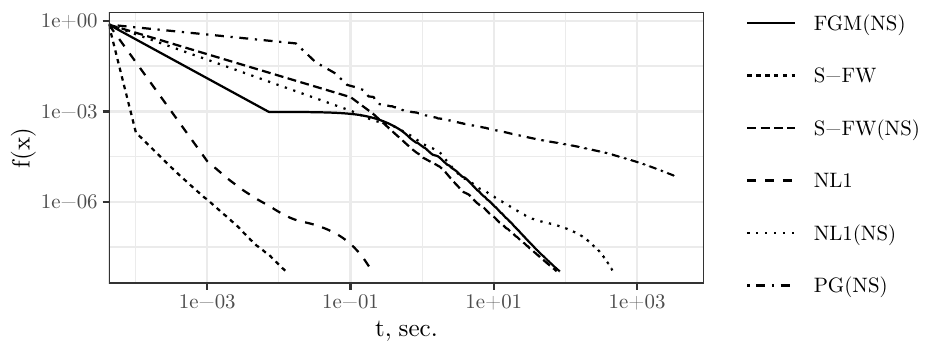}
  \end{center}

  \caption{Time complexity for PageRank over the web-Stanford dataset. The \texttt{S-FW} is fastest one for both sparse and non-sparse (NS) method versions. Sparse algorithms significantly mostly outperforms non-sparse versions.}
  \label{fig:web_stanford}
\end{figure}

\begin{figure}[h!]
  \begin{center}
  \includegraphics[width=\linewidth]{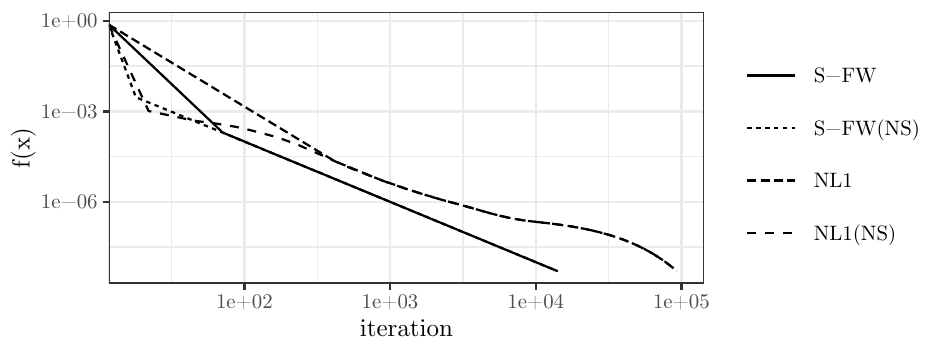}
  \end{center}

  \caption{Iteration complexity for PageRank over the web-Stanford dataset. Sparse methods have almost the same iteration complexity as the non-sparse ones. 
  }
  \label{fig:web_stanford_iter}
\end{figure}

\begin{figure}[h!]
  \begin{center}
  \includegraphics[width=\linewidth]{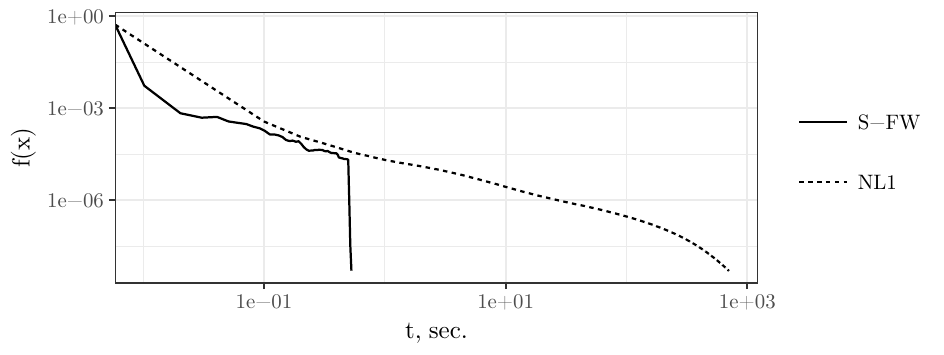}
  \end{center}

  \caption{Time complexity for PageRank over the web-BerkStan dataset. The \texttt{S-FW} algorithm significantly outperforms the \texttt{NL1} algorithm and has a significant gain in the vicinity of the optimal point. }
  \label{fig:web_notredame}
\end{figure}

\begin{figure}[h!]
  \begin{center}
  \includegraphics[width=\linewidth]{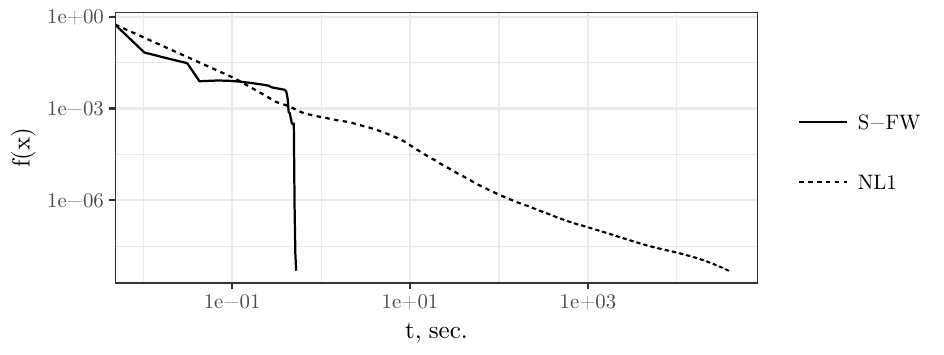}
  \end{center}
  
  \caption{Time complexity for PageRank over the web-BerkStan dataset. The \texttt{S-FW} algorithm significantly outperforms the \texttt{NL1} algorithm and has a significant gain in the vicinity of the optimal point. }
  \label{fig:web_berkstan}
\end{figure}

\begin{figure}[h!]
  \begin{center}
  \includegraphics[width=\linewidth]{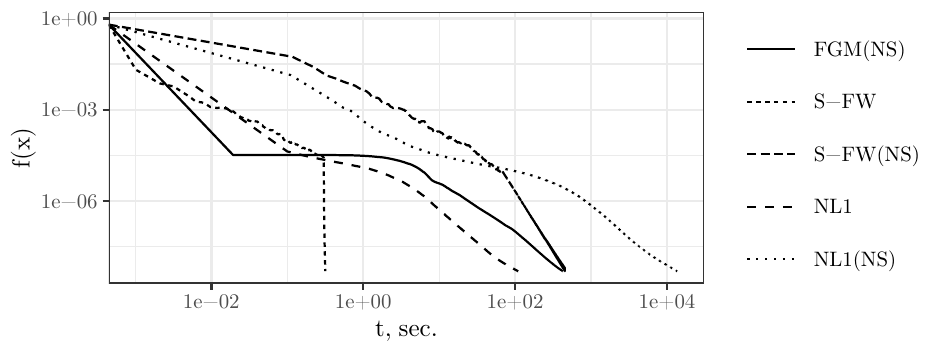}
  \end{center}

  \caption{Time complexity for PageRank over the web-Google dataset. Sparse \texttt{S-FW} method is the fastest among compared. The convergence of  the Similar triangles algorithm (\texttt{FGM-NS})  \cite{nesterov2018lectures,gasnikov2018universal} is faster in terms of the number of iterations, but slower in terms of the computational time.}
  \label{fig:web_google}
\end{figure}

\begin{figure}[h!]
  \begin{center}
  \includegraphics[width=\linewidth]{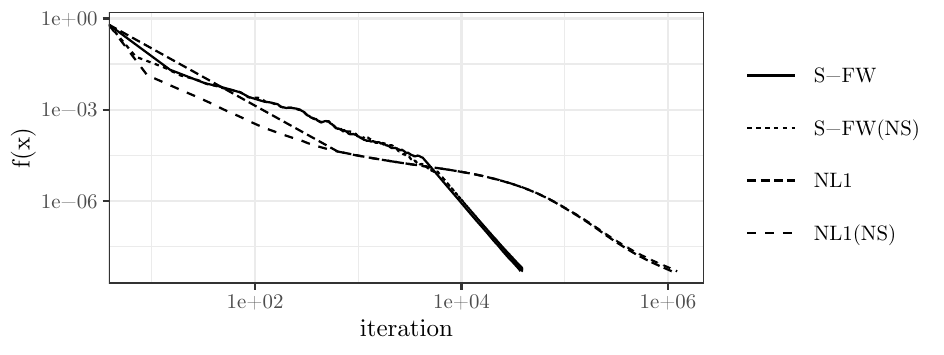}
  \end{center}

  \caption{Iteration complexity for PageRank over the web-Google dataset. The convergence of method's sparse versions is very close to the non-sparse (NS) ones. }
  \label{fig:web_google_iter}
\end{figure}

We also implement and test some other non-sparse (NS) methods and non-sparse versions of our \texttt{S-FW} and \texttt{NL1} methods. Non-sparse versions do not perform any sparse updates and replace them with ``classic'' full-update operations, so we can check not only the speedup of proposed sparse updates but also it's computational accuracy and stability.

We implement a standard non-sparse Projected Gradient Method (PG) \cite{nesterov2018lectures}, and the Similar Triangles algorithm \cite{gasnikov2018universal,nesterov2018lectures}, which we further refer as \texttt{FGM-NS}. 
%
%
Test results prove the accuracy and stability of our sparse methods -- in all cases proposed methods outperforms non-sparse ones and operates very closely it's ``full'' versions (\texttt{S-FW} vs. \texttt{S-FW(NS)} for example). Details are provided in Table~\ref{tbl:a_webgraphs_ns}, 
Figures~\ref{fig:web_stanford}, \ref{fig:web_stanford_iter}, \ref{fig:web_google}, and \ref{fig:web_google_iter}.
\end{enumerate}

\section{Conclusion}
In this paper, we have proposed three novel algorithms to solve the PageRank problem. All the algorithms can be viewed as guided versions of coordinate or block-coordinate descent and demonstrate superior practical performance. In further works, the authors intend to devote more attention to sparsification techniques and interplay between the problem sparsity, dimension, and desired accuracy. 

The work of Anton Anikin was supported by RFBR No. 18-29-03071 mk. The work of Alexander Gasnikov was supported by the Ministry of Science and Higher Education of the Russian Federation (Goszadaniye) No. 075-00337-20-03, project No. 0714-2020-0005. The work of Alexander Gornov was supported by the Ministry of Science and Higher Education of the Russian Federation (Goszadaniye) No. AAAA-A17-117032210080-7. The work of Yury Maximov at LANL was partially funded by the Advanced Grid Modeling Program of the U.S. Department of Energy Office of Electricity (Agreement No. 36306), and LANL-LDRD projects (20210078DR, 20190059DR). The work of Dmitry Kamzolov was funded by RFBR, project No. 19-31-90170. The work of Yurii Nesterov was partially financed by the European Research Council Advanced Grant No. 788368.

\qed

\bibliographystyle{gOMS}
\bibliography{gOMSguide}

\newpage
\appendix
\section{Missing Proofs in Section \ref{sec:NL1}}
Recall the definition of $f_\gamma(x)$:
\begin{gather*}
	f_\gamma(x) = \frac{1}{2} \|Ax\|_2^2 + \sum_{i = 1}^n\frac{\gamma}{2} \left(- x^i
    \right)_+^2. 
\end{gather*}
\begin{lemma}\label{lem:pgm-round}
Let $x_*$ satisfies $e^\top x^* = 1$ and 
$f_{\gamma}(x_*)\le \varepsilon^2$
for some $\gamma > 0$. 
Then for $\hat x = (x_*)_+ / e^\top (x_*)_+$ we have 
\[
\|A\hat x\|_2^2 \le 4 (1+ 4\gamma^{-1}) \varepsilon^2 
\]
and $\hat x\in \Delta^n_1 =\{x \in \mathbb{R}^n:\, \sum_{i=1}^n x_i = 1, \; x_i \ge 0\}$
\end{lemma}
\begin{proof}
By the conditions of the lemma
\begin{gather}\label{eq:aux00}
	\frac{1}{2} \|Ax_*\|_2^2 + \sum_{i = 1}^n\frac{\gamma}{2} \left(-x^i_*
    \right)_+^2 \le \varepsilon^2.
\end{gather}
Let $x = (x_*)_+ - (-x_*)_+$, then by the triangle inequality we have
\begin{gather}\label{eq:aux-1}
	\|(A(x_*)_+)\|_2 \le  \|Ax_*\|_2 + \|(A(-x_*)_+)\|_2 \le \sqrt{2} \varepsilon + \|(A(-x_*)_+)\|_2.
\end{gather}
By the Perron-Frobenius theorem, $|\lambda_i(P^\top)| \le 1$, $1\le i \le n$. Thus 
$\lambda_{\max{}}(A) = \lambda_{\max{}} (I - P^\top) \le 2$.
By Inequality~\ref{eq:aux00} one has
\begin{gather}\label{eq:aux-2}
\|A(-x_*)_+\|_2 \le \|(-x_*)_+\|_2 \le 2 \sqrt{\frac{2}{\gamma}} \varepsilon.
\end{gather}
Using Inequalities~\ref{eq:aux00}, \ref{eq:aux-1} and \ref{eq:aux-2} we have the final estimate
\[
  \|A(x_*)_+\|_2 \le \sqrt{2}\varepsilon + \|A(- x_*)_+\|_2 \le \sqrt{2}\varepsilon + 2\sqrt{\frac{2}{\gamma}} \varepsilon. 
\]
By definition of $\hat x$, we have $\|A\hat x\|_2 ((x_*)_+^\top e) = \|A(x_*)_+\|_2$. Since $e^\top((x_*)_+ - (-x_*)_+) = 1$, $e^\top x_* = 1$, and $(-x_*)_+^\top e \ge 0$ we have 
$(x_*)_+^\top e \ge 1$ and $\|A\hat x\|_2 \le \|A(x_*)_+\|_2$. An application of the Cauchy-Schwartz inequality concludes the proof of the lemma. 
\end{proof}

Proposition B.1 of \cite{allen2014linear} establishes convergence rate of the gradient descent in arbitrary norm. 
\begin{prop}[Proposition B.1 of \cite{allen2014linear}]\label{prop:allen}
Let $f(x)$ be a convex, differentiable function that is $L$-smooth with respect to $\|\cdot\|$ on $Q = \mathbb{R}^n$, and $x_0$ any initial point in $Q$. Consider the sequence of $N$ gradient steps $x_{k+1} = \argmin_{y\in Q} \left\{ \frac{L}{2} \|y-x\|^2 + \nabla f(x_k)^\top (y-x)\right\}$, then the last point $x_{N}$ satisfies 
\[f(x_{N}) - f(x_*) \le 2\frac{LR^2}{N},\]
where $R = \max_{x:\, f(x) \le f(x_0)} \|x-x_*\|$, and $x_*$ is any minimizer of $f$. 
\end{prop}

Now we are ready to proof Theorem~\ref{thm:NL1}. 
\begin{proof}[Proof of the Theorem \ref{thm:NL1}]
A single iteration of Algorithm~\ref{algo:l1-pgm} results in a sparse update vector $h$ containing at most two non-zero coordinates which correspond to the minimal and the maximal coordinates of the gradient (see Lemma~\ref{lem:nl1-min} for the details). Then the gradient update for $x_{k+1} = x_k + h_k$ is
\[
\nabla f_\gamma(x_{k} +  h_k) = \nabla f_\gamma(x_k) + \eta A^\top A h_k + \gamma (- x_k - h_k)_+ - \gamma (- x_k)_+
\]
and requires $O(d^2 \log (n/d^2 + 2))$ arithmetic operations by using a set of $\lceil n/d^2\rceil$ binary heaps described earlier in Section~\ref{sec:NL1}. An update to the  function value is
\[
f_\gamma (x_k + h_k) = f_\gamma (x_k) + 2 h_k^\top A x_k + \|Ah_k\|_2^2/2 + \frac{\gamma}{2} (-x_k - h_k)_+^2 - \frac{\gamma}{2} (-x_k)_+^2
\]
and similarly requires at most $O(d^2 \log (n/d^2 + 2))$ operations as $h_k$ contains at most 2 non-zero coordinates. Notice that the size of the level set 
\[R = \max_{x:\, f(x) \le f(x_0)} \|x-x_*\|_1\]
at a point $x$ is bounded from above as $R \le 4\sqrt{2n f_\gamma(x_0)/\gamma} + 2$. Indeed 
\begin{align}\label{eq:aux-3}
\frac{\gamma}{2} \frac{\|(-x)_+\|_1^2}{n}  \le \frac{\gamma}{2} \|(-x)_+\|_2^2 & = \frac{\gamma}{2} \sum_{i=1}^n(-x^i)_+^2 \nonumber\\ & \le \frac{1}{2}\|Ax\|^2_2 + \frac{\gamma}{2} \sum_{i=1}^n(-x^i)_+^2  = f_\gamma(x)\le f_\gamma(x_0)
\end{align}
and 
\begin{align}\label{eq:aux-4}
R = 
\max_{x:\, f(x) \le f(x_0)} \|x-x_*\|_1 
& \le \max_{x: f_\gamma(x) \le f_\gamma(x_0)} 2\|x\|_1 \nonumber\\
& \le 4 \|(-x)_+\|_1 + 2 \le 4 \sqrt{\frac{2n f_\gamma(x_0)}{\gamma}} + 2,
\end{align}
where $\|x\|_1 \le \|(-x)_+\|_1 + \|(x)_+\|_1 \le 2\|(-x)_+\|_1 + 1$ since $\sum_{i=1}^n x^i = 1$. 

The remainder of the proof will consists of two phases. First, we estimate the time complexity of the algorithm to achieve $f_\gamma(x_{(k)}) \le \gamma/(8n)$. After that, starting from~$x_{(k)}: \, f_\gamma(x_{(k)}) \le \gamma/(8n)$ we find the complexity of the algorithm to achieve~$f_\gamma(x_{m})~\le~\varepsilon^2$. 

Now fix any $\delta > 0$ such that $n \delta^2 < 1$, and let $\varepsilon_k^2 = (\delta^2 n)^k f_\gamma(x_0)$, $k \ge 1$.
To achieve $f_\gamma(x_{(1)}) \le n\delta^2 f_\gamma(x_0)$ 
one needs according to Proposition~\ref{prop:allen} and Eq.~\eqref{eq:aux-4} at most:
\begin{align*}
T_{(1)} = 16\frac{(1+ \gamma)(1 + 8n f_\gamma(x_0)/\gamma)}{\varepsilon_1^2} = 16\frac{(1+\gamma)(1/n + 8f_\gamma(x_0)/\gamma)}{\delta^2} \le 256 \frac{(1+\gamma) f_\gamma(x_0)}{\gamma \delta^2}
\end{align*}
iterations. Similarly, accuracy $f(x_{(k)}) \le \varepsilon^2_k$ can be achieved in 
\begin{align*}
    T_{(k)} & \le 16\frac{(1 + \gamma)(1 + 8n f_\gamma(x_{(k)})/\gamma)}{\varepsilon_k^2}   
    \le 256\frac{(1 + \gamma)n f_\gamma(x_{(k-1)})}{\gamma \varepsilon_k^2} \\
    & = 256 \frac{(1+\gamma) n  (\delta^2 n)^{k-1}f_\gamma(x_0)}{\gamma \delta^{2k}n^{k}} 
    = 256 \frac{(1+\gamma) f_\gamma(x_0)}{\gamma \delta^{2}}  = T_{(1)}
\end{align*}
starting from $x_{(k-1)}$, s.t. $f_\gamma(x_{(k-1)})\le \varepsilon_{k-1}^2 f_\gamma(x_0)$. 

%
%
%

Then $f_\gamma(x_{(k)}) \le \gamma/ (8n)$ in at most $k$ restarts: 
\[8n \gamma f_\gamma(x_0)/\gamma \le 8n (\delta^2 n)^k f_\gamma(x_0)/\gamma~\le~1,\]
e.g. $k = O(\log (n f_\gamma(x_0)/\gamma)/\log (n\delta^2))$. 
Thus the overall number of iterations is 
\[
    T = \sum_{i=1}^k T_{(i)} = O\left(\frac{\log(n f_\gamma(x_0)/\gamma)}{\delta^2 \log (n\delta^2)}\right).
\]

The remaining time required to solve the problem starting with $x_{(k)}$ is bounded from above by Proposition~\ref{prop:allen} as $32 (1+\gamma)/\varepsilon^2$. Thus, the overall complexity of the algorithm does not exceed
\begin{gather}\label{eq:aux_est}
T  = O\left(
    n + d^2 \log (n/d^2 + 2)
    \left(\frac{1+\gamma}{\varepsilon^2} + \inf_{\delta: n\delta^2 < 1} \left[\frac{\log (n f_\gamma(x_0)/\gamma)}{\delta^2 \log (n\delta^2)}\right]\right)
    \frac{(1+\gamma)^2}{\gamma}
    \right).
\end{gather}
by Lemma \ref{lem:pgm-round}. Taking $\gamma = O(1)$, and notice that $f_\gamma(x_0) = 1$ if $x_0$ is a one of the simplex vertices one has 
%
\begin{align*}
T  & = O\left(
    n + d^2 \log (n/d^2 + 2)
    \left(\frac{1}{\varepsilon^2} + \inf_{\delta: n\delta^2 < 1} \left[\frac{\log n}{\delta^2 \log (n\delta^2)}\right]\right)\right) \\
    & = 
    O\left(
    n + \frac{d^2 \log (n/d^2 + 2)}{\varepsilon^2}
    + n d^2 \log (n/d^2 + 2) \log n
    \right),
\end{align*}
\noindent which completes the proof of the theorem. 
\end{proof}

\section{Missing Proofs in Section \ref{sec:GK}}\label{sec:appB}

Theorem~\ref{thm:GK-ch-p} provides us with an upper bound on the efficiency of this strategy. Here we assume that the number of iterations $N$ is known in advance. Let us emphasize that Algorithm~\ref{algo:GK} is a  version of the Mirror Descent algorithm with randomized projecting aims to support sparse updates (see Section~\ref{sec:GK} for details). 

Our main tool below is Proposition~\ref{prop:matrix-game}, which  establishes the convergence rate of stochastic online optimization for linear functions $f_k$ linear in $x$. We refer to recent results in \cite{gasnikov2015efficiency} for a more general problem setup. 

In the proof of the following proposition we mostly follow \cite{cesa2006prediction} and \cite{gasnikov2015efficiency}. Our proof is based on the recent results for 
dual averaging method \cite{nesterov2009primal} that gives essentially the same sequence of steps as the mirror descent for simplex constrained convex optimization problems~\cite{allen2014linear}.  
\begin{prop}\label{prop:matrix-game}
Let $\{f_k(x) = x^\top A y_k\}_{k=1}^N$ be a set of functions $f: \, \mathbb{R}^n \to \mathbb{R}$ of variable $x$ such that
$\|\nabla f_k(x)\|_\infty \le M$ almost surely. Then for a constant step-size policy $\gamma =  M^{-1}\sqrt{2 \log n /N}$ we have
\[\psi_N \doteq \frac{1}{N} \sum_{k=1}^N f_k(x_k) \]
\[
	\mathbb{E}_N
	\left[\psi_N  \right] - \min_{x\in \Delta_1^n}  \frac{1}{N}\sum_{k=1}^N f_k(x) \le M\sqrt{\frac{ 2\log n}{N}}.
\]
where $x_k$ is given as 
\begin{gather}\label{eq:upd3}
x_{k} = e_i, \text{ with probability } p_{k}^i = \frac{p^i_{k-1} \exp\left( - \gamma \frac{\partial f_{k-1} (x_{k-1})}{\partial x^i} \right)}{\sum_{t=1}^n p_{k-1}^t\exp\left( -  \gamma \frac{\partial f_{k-1} (x_{k-1})}{\partial x^t} \right)}, \; k\ge 1
\end{gather}
\noindent and $p_0^i = 1/n$ for all $i$, $1\le i \le n$, and the expectation $\mathbb{E}_N$ is taken over the choice of $x_1, \dots, x_N$. 
%
Moreover, for any $\Omega > 0$
\[
Prob\left[\psi_N > M\sqrt{\frac{2}{N}} \left( \sqrt{\log n} + 2\sqrt{ \Omega}\right) \right] \le \exp (-\Omega). 
\]
\end{prop}
\begin{proof}
Let $\omega(x) = \sum_{i=1}^n x_i \log x_i$ be a distance generating functions, and $\omega^*(z) = \log \left\{\sum_{j=1}^n \exp(z^j)\right\}$ be its conjugate. Then the step of the dual averaging algorithm gives: 
\begin{align}
    \label{eq_zk}
   & {z}_{k} = {z}_{k-1} - \gamma \nabla f_k(x_{k}),\\
   \label{eq_xk}
   &x_{k+1} = \nabla \omega^* ({z}_{k}), 
\end{align}
 
Instead, our update rule uses a randomised projection $x_{k+1}$ of ${z}_{k}$ on a unit simplex, according to Eq.~\eqref{eq:upd3}.

Let ${z}_{k,t} = {z}_{k} t + (1-t) {z}_{k-1}$, $t\in \mathbb{R}$,
then 
\begin{align}\label{eq:aux01}
\omega^*({z}_{k}) & = \omega^*({z}_{k-1}) + \int_0^1 ({z}_{k} - {z}_{k-1})^\top \nabla \omega^* (t {z}_{k}  + (1-t){z}_{k-1}) d t \nonumber\\
& = \omega^*(y_{k-1}) -\gamma \nabla f_k(x_k)^\top \nabla \omega^*({z}_{k-1})  \nonumber \\ 
& \hspace{7mm} -  \gamma \nabla f_k(x_k) \int_0^1  \left[\nabla \omega^* ({z}_{k,t}) - \nabla \omega^*({z}_{k-1}) \right]d t \nonumber\\
& \le \omega^*({z}_{k-1}) - \gamma \nabla f_k(x_k)^\top \nabla \omega^*({z}_{k-1}) \nonumber\\ 
& \hspace{7mm}+ \gamma \|\nabla f_k (x_k)\|_\infty \int_0^1 \left\|\nabla \omega^* ({z}_{k,t}) - \nabla \omega^*({z}_{k-1}) \right\|_1 d t,
\end{align}
where the last is due to Hoelder's inequality. 
By the 1-strong convexity of $\omega(x)$ with respect to the $\ell_1$-norm we have 
\[
    \left\|\nabla \omega^*({z}') - \nabla \omega^*({z})\right\|_1 \le \|{z}' - {z}\|_\infty. 
\]
Then by Inequality~\eqref{eq:aux01} and \eqref{eq_zk} we have
\begin{equation}
\label{eq:aux02}
\begin{aligned}
\omega^*({z}_k) &\le \omega^*({z}_{k-1}) - \gamma \nabla f_k (x_k)^\top \nabla \omega^*({z}_{k-1}) + \frac{\gamma^2 \|\nabla f_k(x_k)\|_\infty^2}{2}\\
&\overset{\eqref{eq_xk}}{\le} \omega^*({z}_{k-1}) - \gamma \nabla f_k (x_k)^\top x_{k}+ \frac{\gamma^2 \|\nabla f_k(x_k)\|_\infty^2}{2}.
\end{aligned}
\end{equation}
Next summing up Ineq.~\eqref{eq:aux02} for all $k$, $1\le k \le {N}$, we have
\begin{gather*}
	\gamma\sum_{k=1}^N  x_k^\top \nabla f_k (x_k) \le -\omega^*(z_N) + \omega^*(z_0) + \frac{\gamma^2}{2}\sum_{k=1}^N  \left\|\nabla f_k(x_k)\right\|_\infty^2.
\end{gather*}
As $z_0 = 0$ and \eqref{eq_zk} we get
\begin{gather*}
	\gamma\sum_{k=1}^N  (x_k-x)^\top \nabla f_k (x_k) \le -\omega^*(z_N) + \omega^*(z_0)+ x^T y_N + \frac{\gamma^2}{2}\sum_{k=1}^N  \left\|\nabla f_k(x_k)\right\|_\infty^2.
\end{gather*}
From $\omega^*(z_0) = 0$ and Young's inequality $\omega^*(z)+\omega(x)\geq x^T z$ we get
\begin{gather*}
	\gamma\sum_{k=1}^N  (x_k-x)^\top \nabla f_k (x_k) \le \omega(x) + \frac{\gamma^2}{2}\sum_{k=1}^N  \left\|\nabla f_k(x_k)\right\|_\infty^2.
\end{gather*}
From convexity of $f$ we get 
\begin{gather*}
	\gamma\sum_{k=1}^N \left(f_k(x_k) - f_k(x)\right) \le \omega(x) + \frac{\gamma^2}{2}\sum_{k=1}^N  \left\|\nabla f_k(x_k)\right\|_\infty^2.
\end{gather*}
Taking expectation $\mathbb{E}_N(\cdot)$ with respect to $x_1, \dots, x_{N}$ (e.g. the choice of $i_1^A, \dots, i_N^A$), we have
\begin{gather*}
	\gamma \sum_{k=1}^N \mathbb{E}_{N}\;\left[f_k(x_k) - f_k(x)\right] \le \omega(x)  + \frac{\gamma^2}{2} \sum_{k=1}^N  \mathbb{E}_{N}\;\left[\left\|\nabla f_k(x)\right\|_\infty^2\right].
\end{gather*}

To finish the proof it remains to note that $\|\nabla f_k(x_k)\|_\infty \le M$ and 
\[
\psi_N \le \min_{\gamma>0}\frac{\log n }{N \gamma}+ \frac{M^2 \gamma}{2}  = M\sqrt{\frac{ 2\log n}{N}}.
\]

The remainder of the proof relies on Azuma's inequality. Let 
\[
Z_j = \sum_{k=1}^j \gamma (x - x_k)^\top \nabla f_k(x_k)
\]
is a Martingale satisfying $|Z_{j+1} - Z_j| \le c_j \doteq 4M \gamma$ almost surely. By Azuma's inequality we have
\[
Prob[Z_N \ge t] \le \exp \left(- \frac{t^2}{2 \sum_{j=1}^N c_j^2} \right). 
\]
Setting $t = 4M\gamma \sqrt{2\Omega N}$ finishes the proof of the proposition. 
\end{proof}

Now we are ready to proof the Theorem~\ref{thm:GK-ch-p}. 
\begin{proof}[Proof of the Theorem~\ref{thm:GK-ch-p}]
For a point $({\bar x}_N, {\bar y}_N)$ defined by Algorithm~\ref{algo:GK} we have:
\begin{align*}
  0 & \le 
  ||A\overline{x}_N||_{\infty}  \overset{\eqref{eq:matrix-game}}{=} \max_{y \in \Delta^{2n}_1} \langle y, \widetilde{A}\overline{x}_N \rangle  \overset{\eqref{eq_saddle_opt}}{=} \max_{y \in \Delta^{2n}_1} \langle y, \widetilde{A}\overline{x}_N \rangle - \max_{y \in \Delta^{2n}_1} \min_{x \in \Delta^{n}_1} \langle y, \widetilde{A}x \rangle \\
   & \le \max_{y \in \Delta^{2n}_1} \langle y, \widetilde{A}\overline{x}_N \rangle - \min_{x \in \Delta^{n}_1} \langle \overline{y}_N, \widetilde{A}x \rangle \\
   & = \left\{\max_{y \in \Delta^{2n}_1} \langle y, \widetilde{A}\overline{x}_N \rangle - \frac{1}{N}\sum\limits_{k=1}^{N} \langle y_k, \widetilde{A}x_k \rangle\right\} + \left\{\frac{1}{N}\sum\limits_{k=1}^{N} \langle y_k, \widetilde{A}x_k \rangle - \min_{x\in \Delta^n_1} \langle \overline{y}_N, \widetilde{A}x \rangle \right\} \\
   & 
  \le\sqrt{\frac{2}{N}}\Big( \sqrt{\ln(2n)} + 2\sqrt{\ln(1/\delta)} \Big) + \sqrt{\frac{2}{N}}\Big( \sqrt{\ln n} + 2\sqrt{\ln(1/\delta)} \Big) \\
  & = \sqrt{\frac{2}{N}}\Big( \sqrt{\ln(2n)}+ \sqrt{\ln(n)} + 4\sqrt{\ln(1/\delta)} \Big),
\end{align*}
where the last estimate is accurate owing to Proposition~\ref{prop:matrix-game} with $M = 1$. That is, with probability at least~$1 - \delta$, it is sufficient to have
\begin{align*}
  N &\ge \frac{2}{\varepsilon^2}\Big( \sqrt{\ln(2n)}+ \sqrt{\ln(n)} + 4\sqrt{\ln(1/\delta)} \Big)^2, \\
  N &\ge \frac{4}{\varepsilon^2}\Big(\ln(2n)+ \ln(n) + 16\ln(1/\delta) \Big)
\end{align*}
iterations of the \texttt{GK} algorithm in order to guarantee $\|A\overline{x}_N\|_{\infty} \le \varepsilon$. Each update of $x$ or $y$ involves an update of no more than $d$ probabilities in vectors $p$ and $\pi$ corresponding to non-zeros in the gradient. Algorithm~\ref{algo:upd} requires $O(\log n)$ time to update each. Therefore, the time-complexity of the algorithm is bounded from above as 
\begin{equation*}
  O\left(n + \frac{d \ln n (\ln n + \ln(\sigma^{-1}))}{\varepsilon^2}\right).
\end{equation*}
\end{proof}

\end{document}